\documentclass[11pt]{article}
\usepackage{}
\usepackage{amssymb}
\usepackage[hang,flushmargin]{footmisc}
\usepackage{multirow}
\usepackage{amsfonts}
\usepackage{amsmath}
\usepackage{graphicx}
\usepackage{amsmath,amssymb,latexsym,color}
\usepackage[mathscr]{eucal}
 \makeatletter
    
    \newcommand{\Rmnum}[1]
    {\expandafter\@slowromancap\romannumeral #1@}
    \makeatother
\def\wz{\tilde}

\newtheorem{thm}{Theorem}[section]

\newtheorem{prop}[thm]{Proposition}
\newtheorem{lemma}[thm]{Lemma}
\newcounter{foo}[subsection]
\newcounter{fooo}[section]
\newtheorem{step}[foo]{Step}
\newtheorem{stepp}[fooo]{Step}
\newtheorem{cor}[thm]{Corollary}

\newtheorem{example}{Example}[section]
\newtheorem{defin}[thm]{Definition}

\newtheorem{remark}{Remark}[section]

\newcommand{\qed}{\hfill\Box\medskip}

\begin{document}

\renewcommand{\baselinestretch}{1.3}
\title{Weakly distance-regular digraphs of valency three, \Rmnum{2}}

\author{
Yuefeng Yang\textsuperscript{a,b}\quad Benjian Lv\textsuperscript{b}\quad
Kaishun Wang\textsuperscript{b}\\
\\{\footnotesize  \textsuperscript{a} \em  School of Science, China University of Geosciences, Beijing, 100083, China}\\{\footnotesize  \textsuperscript{b} \em Sch. Math. Sci. {\rm \&} Lab. Math. Com. Sys., Beijing Normal University, Beijing, 100875, China }  }
\date{}
\maketitle
\footnote{\scriptsize
{\em E-mail address:} yangyf@mail.bnu.edu.cn(Y.Yang), bjlv@bnu.edu.cn(B.Lv), wangks@bnu.edu.cn(K.Wang).}

\begin{abstract}

In this paper, we classify commutative weakly distance-regular digraphs of valency $3$ with girth more than $2$ and one type of arcs. Combining \cite[Theorem 1.2]{HS04}, \cite[Theorem 1.3]{KSW04} and \cite[Theorem 1]{YYF16}, commutative weakly distance-regular digraphs of valency $3$ are completely determined.



\end{abstract}

\section{Introduction}

Throughout this paper $\Gamma$ always denotes a finite simple digraph. We write $V\Gamma$ and $A\Gamma$ for the vertex set and arc set of $\Gamma$, respectively. A \emph{path} of length $r$ from $x$ to $y$ is a sequence of vertices $(x=w_{0},w_{1},\ldots,w_{r}=y)$
such that $(w_{t-1}, w_{t})\in A\Gamma$ for $t=1,2,\ldots,r$. A digraph is said to be \emph{strongly connected} if, for any distinct vertices $x$ and $y$, there is a path from $x$ to $y$. The length of a shortest path from $x$ to $y$ is called the \emph{distance} from $x$ to $y$
in $\Gamma$, denoted by $\partial(x,y)$. Let $\wz{\partial}(x,y)=(\partial(x,y),\partial(y,x))$. We call $\wz{\partial}(x,y)$ the \emph{two way distance} from $x$ to $y$ in $\Gamma$. Denote $\Gamma_{\wz{i}}=\{(x,y)\in V\Gamma\times V\Gamma\mid\wz{\partial}(x,y)=\wz{i}\}$ for any $\wz{i}\in\wz{\partial}(\Gamma)$, where $\wz{\partial}(\Gamma)=\{\wz{\partial}(x,y)\mid x,y\in V\Gamma\}$. An arc $(u,v)$ of $\Gamma$ is of \emph{type} $(1,r)$ if $\partial(v,u)=r$. A path $(w_{0},w_{1},\ldots,w_{r-1})$ is said to be a \emph{circuit} of length $r$ if $\partial(w_{r-1},w_{0})=1$. A circuit of minimal length is said to be a \emph{minimal circuit}. The \emph{girth} of $\Gamma$ is the length of a minimal circuit.

A strongly
connected digraph $\Gamma$ is said to be \emph{weakly distance-regular} if, for any $\wz{h}$, $\wz{i}$, $\wz{j}\in\wz{\partial}(\Gamma)$, the cardinality of set $$P_{\wz{i},\wz{j}}(x,y):=\{z\in V\Gamma\mid\wz{\partial}(x,z)=\wz{i}~\textrm{and}~\wz{\partial}(z,y)=\wz{j}\}$$ is constant whenever $\wz{\partial}(x,y)=\wz{h}$. This constant is denoted by $p_{\wz{i},\wz{j}}^{\wz{h}}$. In other words, $\Gamma$ is weakly distance-regular if $(V\Gamma,R)$ is an association scheme (\cite{EB84,PHZ96,PHZ05}), where $R=\{\Gamma_{\wz{i}}\mid\wz{i}\in\wz{\partial}(\Gamma)\}$. The integers
$p_{\wz{i},\wz{j}}^{\wz{h}}$ are called the \emph{intersection numbers}. We say
that $\Gamma$ is \emph{commutative} if $p_{\wz{i},\wz{j}}^{\wz{h}}=p_{\wz{j},\wz{i}}^{\wz{h}}$ for all $\wz{i}$, $\wz{j}$,
$\wz{h}\in\wz{\partial}(\Gamma)$. We call that $\Gamma$ is \emph{primitive} if the digraph $(V\Gamma,\Gamma_{\wz{i}})$ is strongly connected for all $\wz{i}\in\wz{\partial}(\Gamma)\setminus\{(0,0)\}$. A  weakly distance-regular digraph is \emph{thin} (resp. \emph{quasi-thin}) if the maximum value of its intersection numbers is $1$ (resp. $2$).

Some special families of weakly distance-regular digraphs were classified. See \cite{KSW03,HS04} for valency $2$, \cite{HS04} for thin case and \cite{YYF16+} for quasi-thin case under the assumption of the commutativity. In \cite[Theorem 2.3]{ZL11}, Li et al. gave the classification of weakly distance-regular digraphs with girth more than $2$ whose any path of length $2$ is contained in a minimal circuit. Also in \cite{HS04}, Suzuki proposed the project
to classify weakly distance-regular digraphs of valency $3$. In \cite{KSW04}, Wang determined all such digraphs with girth $2$ under the assumption of the commutativity. In \cite{YYF16}, we classified such digraphs with girth more than $2$ and two types of arcs. In this paper, we continue this project, and obtain the following result.

\begin{thm}\label{Main}
Let $\Gamma$ be a commutative weakly distance-regular digraph of valency $3$ and girth $g$, where $g\geq3$. If $\Gamma$ has one type of arcs, then $\Gamma$ is
isomorphic to one of the following digraphs:\vspace{-0.3cm}
\begin{itemize}
\item [{\rm(i)}] ${\rm Cay}(\mathbb{Z}_{7},\{1,2,4\})$.\vspace{-0.3cm}

\item [{\rm(ii)}] ${\rm Cay}(Q_{8},\{i,j,k\})$, $Q_{8}$ is the quaternion group of order $8$.\vspace{-0.3cm}

\item [{\rm(iii)}] ${\rm Cay}(\mathbb{Z}_{13},\{1,3,9\})$.\vspace{-0.3cm}

\item [{\rm(iv)}] The eighteenth digraph with $18$ vertices in {\rm\cite{H}}.\vspace{-0.3cm}

\item [{\rm(v)}] ${\rm Cay}(\mathbb{Z}_{3}\times\mathbb{Z}_{3}\times\mathbb{Z}_{3},\{(1,0,0),(0,1,0),(0,0,1)\})$.\vspace{-0.3cm}

\item [{\rm(vi)}] ${\rm Cay}(\mathbb{Z}_{g}\times\mathbb{Z}_{3},\{(1,0),(1,1),(1,2)\})$.\vspace{-0.3cm}

\item [{\rm(vii)}] ${\rm Cay}(\mathbb{Z}_{n}\times\mathbb{Z}_{n},\{(1,0),(0,1),(n-1,n-1)\})$, $n\notin 3\mathbb{Z}\setminus\{3\}$.\vspace{-0.3cm}

\item [{\rm(viii)}] ${\rm Cay}(\mathbb{Z}_{n}\times\mathbb{Z}_{3n},\{(0,1),(1,1),(n-1,3n-2)\})$, $n\geq2$.
\end{itemize}
\end{thm}

Routinely, all digraphs in above theorem are commutative weakly distance-regular. The digraph in (iv) is not a Cayley digraph by \cite{AH}. For the last three families of Cayley digraphs, in Table 1, we list the two way distance from the identity element to any other element of the corresponding group.

In \cite{MH99,MH00}, Hirasaka gave the classification of primitive commutative association schemes with a non-symmetric relation of valency $3$ or $4$. \cite[Proposition 1.4 and Theorem 1.5]{MH00} and Theorem \ref{Main} imply the following immediately.

\begin{cor}\label{primitive}
If $\Gamma$ is a primitive commutative weakly distance-regular digraph of valency $3$ and girth more than $2$, then $\Gamma$ is
isomorphic to one of the digraphs {\rm(i), (iii)} and {\rm(vii)} in Theorem {\rm\ref{Main}}.
\end{cor}

In order to prove our main theorem, we prepare some basic notations and lemmas.

Let $\Gamma$ be a weakly distance-regular digraph. For each $(i,j)\in\wz{\partial}(\Gamma)$, let $A_{i,j}$ denote a matrix with rows and columns indexed by $V\Gamma$ such that $(A_{i,j})_{x,y}=1$ if $\wz{\partial}(x,y)=(i,j)$, and $(A_{i,j})_{x,y}=0$ otherwise. It follows from the definition of association schemes that
\begin{eqnarray}
A_{i,i'}A_{j,j'}=\sum_{(h,h')\in\wz{\partial}(\Gamma)}p_{(i,i'),(j,j')}^{(h,h')}A_{h,h'}.\nonumber
\end{eqnarray}
For two nonempty subsets $E$ and $F$ of $R$, define
\begin{eqnarray}
EF=\{\Gamma_{\wz{h}}\mid\sum_{\Gamma_{\wz{i}}\in E}\sum_{\Gamma_{\wz{j}}\in F}p_{\wz{i},\wz{j}}^{\wz{h}}\neq0\}, \nonumber
\end{eqnarray}
and write $\Gamma_{\wz{i}}\Gamma_{\wz{j}}$ instead of $\{\Gamma_{\wz{i}}\}\{\Gamma_{\wz{j}}\}$. The size of $\Gamma_{\wz{i}}(x):=\{y\in V\Gamma\mid\wz{\partial}(x,y)=\wz{i}\}$ depends only on $\wz{i}$, denoted by $k_{\wz{i}}$. For any $(a,b)\in\wz{\partial}(\Gamma)$, we usually write $k_{a,b}$ (resp. $\Gamma_{a,b}$) instead of $k_{(a,b)}$ (resp. $\Gamma_{(a,b)}$). Now we list basic properties of intersection numbers which are used frequently in this paper.

\begin{lemma}\label{jiben}
{\rm (\cite[Chapter \Rmnum{2}, Proposition 2.2]{EB84} and \cite[Proposition 5.1]{ZA99})} For each $\wz{i}:=(a,b)\in\wz{\partial}(\Gamma)$, define $\wz{i}^{*}=(b,a)$. The following hold:\vspace{-0.3cm}
\begin{itemize}
\item [{\rm(i)}] $k_{\wz{d}}k_{\wz{e}}=\sum_{\wz{f}\in\wz{\partial}(\Gamma)}p_{\wz{d},\wz{e}}^{\wz{f}}k_{\wz{f}}$.\vspace{-0.3cm}

\item [{\rm(ii)}] $p_{\wz{d},\wz{e}}^{\wz{f}}k_{\wz{f}}=p_{\wz{f},\wz{e}^{*}}^{\wz{d}}k_{\wz{d}}=p_{\wz{d}^{*},\wz{f}}^{\wz{e}}k_{\wz{e}}$.\vspace{-0.3cm}

\item [{\rm(iii)}] $|\Gamma_{\wz{d}}\Gamma_{\wz{e}}|\leq{\rm gcd}(k_{\wz{d}},k_{\wz{e}})$.\vspace{-0.3cm}

\item [{\rm(iv)}] $\sum_{\wz{e}\in\wz{\partial}(\Gamma)}p_{\wz{d},\wz{e}}^{\wz{f}}=k_{\wz{d}}$.\vspace{-0.3cm}

\item [{\rm(v)}] $\sum_{\wz{f}\in\wz{\partial}(\Gamma)}p_{\wz{d},\wz{e}}^{\wz{f}}p_{\wz{g},\wz{f}}^{\wz{h}}=\sum_{\wz{l}\in\wz{\partial}(\Gamma)}p_{\wz{g},\wz{d}}^{\wz{l}}p_{\wz{l},\wz{e}}^{\wz{h}}$.\vspace{-0.3cm}

\item [{\rm(vi)}] ${\rm lcm}(k_{\wz{d}},k_{\wz{e}})\mid p_{\wz{d},\wz{e}}^{\wz{f}}k_{\wz{f}}$.
\end{itemize}
\end{lemma}

Let $\Gamma$ be a digraph satisfying the conditions in Theorem \ref{Main} and $s=\max\{j\mid p_{(1,g-1),(1,g-1)}^{(2,j)}\neq0\}$ throughout the remaining of this paper. Since $p_{(1,g-1),(1,g-1)}^{(2,g-2)}\neq0$, $s\geq g-2$.

A major tool in the proof of our main result is the $h$-chain, which is inspired by the chain defined by Hirasaka in \cite{MH99,MH00}. In $\Gamma$, we call that a path $L:=(y_{0},y_{1},\ldots,y_{j})$ is an \emph{$h$-chain} if $\wz{\partial}(y_{i},y_{i+2})=(2,h)$ for $0\leq i\leq j-2$. If $L$ is a shortest $h$-chain such that $y_{0}=y_{j}$, we say that $L$ is an {\em$h$-line}.

For $\wz{f}\in\wz{\partial}(\Gamma)$, we
define a relation $\wz{\Gamma}_{\wz{f}}$ on the set of $h$-chains as follows. For any two $h$-chains $L=(y_{0},y_{1},\ldots,y_{j})$ and $L'=(y_{0}',y_{1}',\ldots,y_{j'}')$, $(L,L')\in\wz{\Gamma}_{\wz{f}}$ if and only if $j=j'$ and
$\wz{\partial}(y_{i},y_{i}')=\wz{f}$ for $0\leq i\leq j$. Let $L_{0},L_{1},\ldots,L_{m}$ be $h$-lines. We say that $(L_{0},L_{1},\ldots,L_{m})$ is a {\em chain of $h$-lines} if $(L_{i},L_{i+1})\in\wz{\Gamma}_{1,g-1}$ for $0\leq i\leq m-1$, and $(L_{i'},L_{i'+2})\in\wz{\Gamma}_{2,h}$
for $0\leq i'\leq m-2$. In particular, we say that a chain of $h$-lines
$(L_{0},L_{1},\ldots,L_{m})$ is an {\em $h$-plane} if $L_{0}=L_{m}$.

\vspace{3ex}

\noindent\textbf{Outline of the proof of Theorem \ref{Main}.} By \cite[Theorem 2.3]{ZL11} and Lemma \ref{jiben} (ii), if $p_{(1,g-1),(1,g-1)}^{(2,g-2)}=3$, then $\Gamma$ is isomorphic to one of the digraphs in (vi). We only need to consider the case that $p_{(1,g-1),(1,g-1)}^{(2,g-2)}=2$ or $1$.

In Section 2, we prove our main result under the assumption that $p_{(1,g-1),(1,g-1)}^{(2,g-2)}=2$. First, we determine all the possible decompositions of $(A_{1,g-1})^{2}$ and $A_{1,g-1}A_{g-1,1}$. Based on these decompositions, we divide our proof into four cases. For the first three cases, we compute the number of vertices of $\Gamma$, and determine these digraphs according to \cite{H}. For the last case, we give some properties of an $s$-chain, and use $s$-chains to construct $\Gamma$.

In Section 3, we prove our main result under the assumption that $p_{(1,g-1),(1,g-1)}^{(2,g-2)}=1$. We begin with determining all the possible decompositions of $(A_{1,g-1})^{2}$.  Based on these decompositions, we divide our proof into three subsections. For the first two subsections, we compute the number of vertices of $\Gamma$, and determine these digraphs according to \cite{H}. For the last subsection, we use $(g-2)$-chains to construct of $\Gamma$.

In Section 4, we give a proof of Lemma \ref{pre-main}.

\section{$p_{(1,g-1),(1,g-1)}^{(2,g-2)}=2$}

In this section, we prove Theorem \ref{Main} under the assumption that $p_{(1,g-1),(1,g-1)}^{(2,g-2)}=2$. We begin with the following lemma.

\begin{lemma}\label{jiben2}
We have $g=3$. Moreover, one of the following holds:\vspace{-0.3cm}
\begin{itemize}
\item [{\rm C1)}] $(A_{1,2})^{2}=2A_{2,1}+A_{1,2}$ and $A_{1,2}A_{2,1}=3I+A_{1,2}+A_{2,1}$.\vspace{-0.3cm}

\item [{\rm C2)}] $(A_{1,2})^{2}=2A_{2,1}+A_{2,3}$ and $A_{1,2}A_{2,1}=3I+A_{2,3}+A_{3,2}$.\vspace{-0.3cm}

\item [{\rm C3)}] $(A_{1,2})^{2}=2A_{2,1}+3A_{2,s}$ and $A_{1,2}A_{2,1}=3I+2A_{3,3}$.\vspace{-0.3cm}

\item [{\rm C4)}] $(A_{1,2})^{2}=2A_{2,1}+A_{2,s}$ and $A_{1,2}A_{2,1}=3I+A_{3,3}$.
\end{itemize}
\end{lemma}
\textbf{Proof.}~For fixed $x_{0,0}\in V\Gamma$, let $(x_{0,0}=x_{g,0},x_{1,0},\ldots,x_{g-1,0})$ be a minimal circuit and $\Gamma_{1,g-1}(x_{0,0})=\{x_{1,0},x_{1,1},x_{1,2}\}$. Since $p_{(1,g-1),(1,g-1)}^{(2,g-2)}=2$, from Lemma \ref{jiben} (ii), we have $p_{(2,g-2),(g-1,1)}^{(1,g-1)}=2$ and $k_{2,g-2}=3$. Without loss of generality, we may assume that $\Gamma_{2,g-2}(x_{0,0})=\{x_{2,0},x_{2,1},x_{2,2}\}$ and $x_{2,i},x_{2,i+1}\in\Gamma_{1,g-1}(x_{1,i})$, where $x_{2,3}=x_{2,0}$ and $i=0,1,2$. Since $x_{1,0}\in\Gamma_{g-2,2}(x_{3,0})$, there exists a vertex $x_{0}\in P_{(1,g-1),(1,g-1)}(x_{1,0},x_{3,0})\setminus\{x_{2,0}\}$. By $x_{0}\in\Gamma_{2,g-2}(x_{0,0})=\{x_{2,0},x_{2,1},x_{2,2}\}$, we obtain $x_{0}=x_{2,1}$. Since $p_{(2,g-2),(g-1,1)}^{(1,g-1)}=2$, we get $(x_{g-1,0},x_{1,1})\notin\Gamma_{2,g-2}$, which implies $x_{g-1,0}=x_{2,0}$. Thus, $g=3$.

We claim that $p_{(1,2),(2,1)}^{(2,2)}=0$. Suppose not. Let $x,y,z$ be three vertices such that $\partial(x,y)=\partial(z,y)=1$ and $\wz{\partial}(x,z)=(2,2)$. By Lemma \ref{jiben} (i) and (vi), we have $(\Gamma_{1,2})^{2}=\{\Gamma_{2,1},\Gamma_{2,2}\}$, which implies that there exist two vertices $w\in P_{(1,2),(1,2)}(z,x)$ and $w'\in P_{(1,2),(1,2)}(x,z)$. Since $p_{(2,1),(2,1)}^{(1,2)}=p_{(1,2),(1,2)}^{(2,1)}=2$, one gets $w'$ or $y\in \Gamma_{2,1}(w)$. By $z\in P_{(1,2),(1,2)}(w',w)$ or $y\in P_{(1,2),(1,2)}(z,w)$, we obtain $\Gamma_{1,2}\in(\Gamma_{1,2})^{2}$, a contradiction. Thus, our claim is valid.

In view of Lemma \ref{jiben} (iii), we have $|\Gamma_{1,2}\Gamma_{2,1}|=3$ or $2$.

\textbf{Case 1.} $|\Gamma_{1,2}\Gamma_{2,1}|=3$.

By the claim, one has $A_{1,2}A_{2,1}=3I+p_{(1,2),(2,1)}^{(h,l)}A_{h,l}+p_{(1,2),(2,1)}^{(l,h)}A_{l,h}$ with $(h,l)\in\{(1,2),(2,3)\}$. In view of Lemma \ref{jiben} (i) and (vi), we get $(A_{1,2})^{2}=2A_{2,1}+p_{(1,2),(1,2)}^{(h,l)}A_{h,l}$ and $p_{(1,2),(1,2)}^{(h,l)}k_{h,l}=p_{(1,2),(2,1)}^{(h,l)}k_{h,l}=3$. Lemma \ref{jiben} (iv) implies that $p_{(1,2),(2,1)}^{(h,l)}=p_{(1,2),(1,2)}^{(h,l)}=1$. Thus, C1 or C2 holds.

\textbf{Case 2.} $|\Gamma_{1,2}\Gamma_{2,1}|=2$.

In view of the claim, one has $A_{1,2}A_{2,1}=3I+p_{(1,2),(2,1)}^{(3,3)}A_{3,3}$. Since $\Gamma_{1,2}\notin(\Gamma_{1,2})^{2}$, by Lemma \ref{jiben} (i) and (vi), we obtains $(A_{1,2})^{2}=2A_{2,1}+p_{(1,2),(1,2)}^{(2,s)}A_{2,s}$ with $s>1$. It follows from Lemma \ref{jiben} (v) and the commutativity of $\Gamma$ that $p_{(1,2),(1,2)}^{(2,1)}p_{(2,1),(2,1)}^{(1,2)}+p_{(1,2),(1,2)}^{(2,s)}p_{(2,1),(2,s)}^{(1,2)}=3+p_{(2,1),(1,2)}^{(3,3)}p_{(3,3),(1,2)}^{(1,2)}$. In view of Lemma \ref{jiben} (i), we get $p_{(1,2),(2,1)}^{(3,3)}k_{3,3}=6$ and $p_{(1,2),(1,2)}^{(2,s)}k_{2,s}=3$, which imply $p_{(3,3),(1,2)}^{(1,2)}=p_{(2,1),(2,1)}^{(1,2)}=2$ and $p_{(2,1),(2,s)}^{(1,2)}=1$ from Lemma \ref{jiben} (ii). Thus, $2p_{(2,1),(1,2)}^{(3,3)}=p_{(1,2),(1,2)}^{(2,s)}+1$. If $p_{(1,2),(1,2)}^{(2,s)}=3$, then $p_{(2,1),(1,2)}^{(3,3)}=2$ and C3 holds; if $p_{(1,2),(1,2)}^{(2,s)}=1$, then $p_{(2,1),(1,2)}^{(3,3)}=1$ and C4 holds.$\qed$

In the following, we divide the proof into four subsections according to separate assumptions based on Lemma \ref{jiben2}.

\subsection{The case C1}

\begin{prop}\label{jcs2h=3(1)}
In this case, $\Gamma$ is isomorphic to the digraph in Theorem {\rm\ref{Main} (i)}.
\end{prop}
\textbf{Proof.}~Since $\Gamma_{1,2}\Gamma_{2,1}=\Gamma_{2,1}\Gamma_{1,2}$ from the commutativity of $\Gamma$, we have $(\Gamma_{1,2})^{i}=\{\Gamma_{0,0},\Gamma_{1,2},\Gamma_{2,1}\}$ for $i\geq3$. Hence, $\wz{\partial}(\Gamma)=\{(0,0),(1,2),(2,1)\}$ and $|V\Gamma|=7$. By \cite{H}, $\Gamma\simeq\textrm{Cay}(\mathbb{Z}_{7},\{1,2,4\})$.$\qed$

\subsection{The case C2}

\begin{prop}\label{jcs2h=3}
In this case, $\Gamma$ is isomorphic to the digraph in Theorem {\rm\ref{Main} (iii)}.
\end{prop}
\textbf{Proof.}~Pick a path $(x,y,z)$ such that $\wz{\partial}(x,z)=(2,3)$. Suppose that $\Gamma_{1,2}(z)=\{w,w_{1},w_{2}\}$. In view of $p_{(2,1),(2,1)}^{(1,2)}=p_{(1,2),(1,2)}^{(2,1)}=2$, we may assume $\wz{\partial}(y,w_{1})=\wz{\partial}(y,w_{2})=(2,1)$. Since $p_{(1,2),(2,1)}^{(2,3)}=1$ and $k_{2,3}=3$ from Lemma \ref{jiben} (i), we get $p_{(2,3),(1,2)}^{(1,2)}=p_{(3,2),(1,2)}^{(1,2)}=1$ by Lemma \ref{jiben} (ii). Without loss of generality, we may assume that $\wz{\partial}(x,w_{1})=\wz{\partial}(w_{2},x)=(3,2)$. Hence, $w\in P_{(1,2),(2,1)}(x,z)$ and $\Gamma_{2,3}\Gamma_{1,2}=\{\Gamma_{2,3},\Gamma_{3,2},\Gamma_{1,2}\}$. By $\Gamma_{1,2}\Gamma_{2,1}=\Gamma_{2,1}\Gamma_{1,2}$ from the commutativity of $\Gamma$, one gets $(\Gamma_{1,2})^{3}=\{\Gamma_{2,3},\Gamma_{3,2},\Gamma_{0,0},\Gamma_{1,2}\}$.

Since $p_{(1,2),(1,2)}^{(2,3)}=1$, we may assume that $\Gamma_{1,2}(w_{1})=\{y,y_{1},y_{2}\}$ and $y_{1}\in P_{(1,2),(1,2)}(w_{1},x)$. The fact that $p_{(3,2),(2,1)}^{(3,2)}=p_{(2,3),(1,2)}^{(2,3)}\neq0$ implies that $y_{2}\in P_{(3,2),(2,1)}(x,w_{1})$. Hence, $\Gamma_{3,2}\Gamma_{1,2}=\{\Gamma_{3,2},\Gamma_{2,1},\Gamma_{1,2}\}$. It follows that $(\Gamma_{1,2})^{4}=\{\Gamma_{2,3},\Gamma_{3,2},\Gamma_{1,2},\Gamma_{2,1}\}$ and $(\Gamma_{1,2})^{i}=\{\Gamma_{0,0},\Gamma_{1,2},\Gamma_{2,1},\Gamma_{2,3},\Gamma_{3,2}\}$ for $i\geq5$. Then $\wz{\partial}(\Gamma)=\{(0,0),(1,2),(2,1),(2,3),(3,2)\}$ and $|V\Gamma|=13$. By \cite{H}, $\Gamma\simeq\textrm{Cay}(\mathbb{Z}_{13},\{1,3,9\})$.$\qed$

\subsection{The case C3}

\begin{prop}\label{jcs2p=3}
In this case, $\Gamma\simeq{\rm Cay}(\mathbb{Z}_{2}\times\mathbb{Z}_{6},\{(0,1),(1,1),(1,4)\})$.
\end{prop}
\textbf{Proof.}~By Lemma \ref{jiben} (i), $k_{2,s}=1$. Pick a path $(x,y,z)$ such that $\wz{\partial}(x,z)=(2,s)$. It follows from Lemma \ref{jiben} (iii) that $|\Gamma_{2,s}\Gamma_{1,2}|=1$. Since $p_{(2,1),(2,1)}^{(1,2)}=p_{(1,2),(1,2)}^{(2,1)}\neq0$, there exists a vertex $w\in P_{(2,1),(2,1)}(y,z)$ such that $\partial(w,x)=s-1$. By $A_{1,2}A_{2,1}=3A_{0,0}+2A_{3,3}$ and $y\in P_{(1,2),(2,1)}(x,w)$, we have $\wz{\partial}(x,w)=(3,3)$ and $s=4$. Hence, $\Gamma_{2,4}\Gamma_{1,2}=\{\Gamma_{3,3}\}$. Since $\Gamma_{1,2}\Gamma_{2,1}=\Gamma_{2,1}\Gamma_{1,2}$ from the commutativity of $\Gamma$, we get $(\Gamma_{1,2})^{3}=\{\Gamma_{0,0},\Gamma_{3,3}\}$.

Let $\Gamma_{1,2}(w)=\{y,y_{1},y_{2}\}$. Since $\wz{\partial}(x,w)=(3,3)$ and $p_{(2,1),(1,2)}^{(3,3)}=p_{(1,2),(2,1)}^{(3,3)}=2$, there exist two vertices $z_{1},z_{2}$ such that $\{z_{1},z_{2}\}=P_{(2,1),(1,2)}(x,w)$. For each $i=1,2$, we have $\wz{\partial}(z_{j_{i}},y_{i})=(2,1)$ for some $j_{i}\in\{1,2\}$ by $p_{(2,1),(2,1)}^{(1,2)}=p_{(1,2),(1,2)}^{(2,1)}=2$. In view of $p_{(1,2),(1,2)}^{(1,2)}=0$, we obtain $\partial(y_{i},x)=2$. Since $\wz{\partial}(x,w)=(3,3)$ and $p_{(1,2),(2,1)}^{(3,3)}=2$, one gets $\wz{\partial}(x,y_{h})=(4,2)$ and $\wz{\partial}(x,y_{l})=(1,2)$ with $\{h,l\}=\{1,2\}$, which imply $\Gamma_{3,3}\Gamma_{1,2}=\{\Gamma_{1,2},\Gamma_{4,2}\}$. Then $(\Gamma_{1,2})^{4}=\{\Gamma_{1,2},\Gamma_{4,2}\}$.

Since $k_{4,2}=1$, from Lemma \ref{jiben} (iii), we get $\Gamma_{4,2}\Gamma_{1,2}=\{\Gamma_{2,1}\}$, which implies $(\Gamma_{1,2})^{5}=\{\Gamma_{2,1},\Gamma_{2,4}\}$. Hence, $(\Gamma_{1,2})^{j+3i}=(\Gamma_{1,2})^{j}$ for $j=2,3,4$ and $i\geq0$. Then $\wz{\partial}(\Gamma)=\{(0,0),(1,2),(2,1),(2,4),(3,3),(4,2)\}$. By Lemma \ref{jiben} (i), one has $k_{3,3}=3$ and $|V\Gamma|=12$. The desired result follows from \cite{H}.$\qed$

\subsection{The case C4}

In this case, by Lemma \ref{jiben} (i), we have $k_{2,s}=3$ and $k_{3,3}=6$. For $j\geq 2$, we set $\Gamma_{1}=\Gamma_{1,2}$ and
\begin{eqnarray}
\Gamma_{j}=\{(x,y)\mid \textrm{there exists an $s$-chain from $x$ to $y$ of length $j$} \}.\nonumber
\end{eqnarray}
Denote $\Gamma_{i}(x)=\{y\mid(x,y)\in\Gamma_{i}\}$ and $k_{i}(x)=|\Gamma_{i}(x)|$ for $i\geq1$ and $x\in V\Gamma$. Since $p_{(2,s),(2,1)}^{(1,2)}=1$ from Lemma \ref{jiben} (ii), one gets $k_{i}(x)\leq3$. Let $k_{i}=k_{i}(x)$, when $k_{i}(x)$ depends only on $i$ and does not depend on the choice of $x$. Note that $\Gamma_{2}=\Gamma_{2,s}$.

\subsubsection{Properties of $s$-chains}

\begin{lemma}\label{jcs2_3neq(3,3)}
Suppose that $(x_{0},x_{1},x_{2},x_{3})$ is an $s$-chain. Then $(x_{0},x_{3})\notin\Gamma_{3,3}$, $\Gamma_{3}=\Gamma_{\wz{\partial}(x_{0},x_{3})}$ and $k_{3}\in\{1,3\}$.
\end{lemma}
\textbf{Proof.}~Let $\Gamma_{1,2}(x_{2})=\{x_{3},x_{3}',x_{3}''\}$. Since $p_{(2,1),(2,1)}^{(1,2)}=p_{(1,2),(1,2)}^{(2,1)}=2$, we obtain  $x_{3}',x_{3}''\in\Gamma_{2,1}(x_{1})$, which implies $x_{3}',x_{3}''\in\Gamma_{3,3}(x_{0})$ by $\Gamma_{1,2}\Gamma_{2,1}=\{\Gamma_{0,0},\Gamma_{3,3}\}$. If $x_{3}\in\Gamma_{3,3}(x_{0})$, from Lemma \ref{jiben} (i), then $p_{(2,s),(1,2)}^{(3,3)}=k_{2,s}k_{1,2}/k_{3,3}=3/2$, a contradiction. Hence, $(x_{0},x_{3})\notin\Gamma_{3,3}$ and $\Gamma_{2,s}\Gamma_{1,2}=\{\Gamma_{3,3},\Gamma_{\wz{a}_{3}}\}$, where $\wz{\partial}(x_{0},x_{3})=\wz{a}_{3}$.

For any $s$-chain $(y_{0},y_{1},y_{2},y_{3})$, similarly, $\wz{\partial}(y_{0},y_{3})\neq(3,3)$. Since $\Gamma_{2,s}\Gamma_{1,2}=\{\Gamma_{3,3},\Gamma_{\wz{a}_{3}}\}$, we have $\wz{\partial}(y_{0},y_{3})=\wz{\partial}(x_{0},x_{3})$ and $\Gamma_{3}\subseteq\Gamma_{\wz{a}_{3}}$. Conversely, for any $(z_{0},z_{3})\in\Gamma_{\wz{a}_{3}}$, there exist two vertices $z_{1},z_{2}$ such that $z_{1}\in P_{(1,2),(2,s)}(z_{0},z_{3})$ and $z_{2}\in P_{(1,2),(1,2)}(z_{1},z_{3})$ from $x_{1}\in P_{(1,2),(2,s)}(x_{0},x_{3})$. Since $\Gamma_{1,2}\Gamma_{2,1}=\{\Gamma_{0,0},\Gamma_{3,3}\}$, we get $\wz{\partial}(z_{0},z_{2})\neq(2,1)$, which implies $\wz{\partial}(z_{0},z_{2})=(2,s)$ and $(z_{0},z_{3})\in\Gamma_{3}$. Then $\Gamma_{3}\subseteq\Gamma_{\wz{a}_{3}}$, so $\Gamma_{3}=\Gamma_{\wz{a}_{3}}$. By $\Gamma_{2,s}\Gamma_{1,2}=\{\Gamma_{3,3},\Gamma_{\wz{a}_{3}}\}$ and Lemma \ref{jiben} (i),(vi), one has $k_{3}=k_{\wz{a}_{3}}\in\{1,3\}$.$\qed$

\begin{lemma}\label{jcs2gamma4}
Let $(x_{0},x_{1},x_{2},x_{3},x_{4})$ be an $s$-chain. If $p_{(1,2),(3,3)}^{\wz{\partial}(x_{0},x_{4})}\neq0$, then $k_{3}=1$.
\end{lemma}
\textbf{Proof.}~Pick a vertex $z_{1}\in P_{(1,2),(3,3)}(x_{0},x_{4})$. Since $\Gamma_{3,3}\in\Gamma_{1,2}\Gamma_{2,1}$, from the commutativity of $\Gamma$, there exists a vertex $z_{3}\in P_{(2,1),(1,2)}(z_{1},x_{4})$. Lemma $\ref{jcs2_3neq(3,3)}$ implies $z_{1}\neq x_{1}$. By $p_{(1,2),(1,2)}^{(2,1)}=p_{(2,1),(2,1)}^{(1,2)}=2$, we may assume $z_{2}\in P_{(1,2),(1,2)}(z_{1},z_{3})\cap P_{(2,1),(2,1)}(z_{3},x_{4})$ and $P_{(1,2),(1,2)}(z_{2},x_{4})=\{z_{3},w_{3}\}$. Since $z_{3}\in P_{(2,1),(1,2)}(z_{1},x_{4})$ and $p_{(2,1),(1,2)}^{(3,3)}=p_{(1,2),(2,1)}^{(3,3)}=1$, we get $\wz{\partial}(z_{1},w_{3})=(2,s)$.

Suppose $k_{3}\neq1$. By Lemma \ref{jcs2_3neq(3,3)}, $k_{3}=3$. In view of Lemma \ref{jiben} (ii), we have
$p_{(3,3),(1,2)}^{\wz{a}_{4}}k_{\wz{a}_{4}}=p_{\wz{a}_{4},(2,1)}^{(3,3)}k_{3,3}$ and $p_{\wz{a}_{3},(1,2)}^{\wz{a}_{4}}k_{\wz{a}_{4}}=p_{\wz{a}_{4},(2,1)}^{\wz{a}_{3}}k_{\wz{a}_{3}}$, where $\wz{\partial}(x_{0},x_{3})=\wz{a}_{3}$ and $\wz{\partial}(x_{0},x_{4})=\wz{a}_{4}$. Since $k_{3,3}=6$ and $k_{\wz{a}_{3}}=k_{3}=3$, one gets $p_{(3,3),(1,2)}^{\wz{a}_{4}}\diagup p_{\wz{a}_{3},(1,2)}^{\wz{a}_{4}}\geq\frac{2}{3}$, which implies $p_{\wz{a}_{3},(1,2)}^{\wz{a}_{4}}=1$ by Lemma \ref{jiben} (iv). In view of $x_{1}\neq z_{1}$ and Lemma \ref{jcs2_3neq(3,3)}, we obtain $(x_{0},w_{3})\notin\Gamma_{3}$. Hence, $\wz{\partial}(x_{0},z_{2})=(2,1)$.

Observe that $(x_{0},x_{4})\in\Gamma_{4}$ and $\Gamma_{2,1}\Gamma_{1,2}=\Gamma_{1,2}\Gamma_{2,1}=\{\Gamma_{0,0},\Gamma_{3,3}\}$. In view of Lemma \ref{jcs2_3neq(3,3)}, one obtains $\Gamma_{2,1}\cap\{(z_{2},x_{1}),(x_{3},z_{2}),(x_{4},x_{0})\}=\emptyset$, which implies that $\{(z_{2},x_{1}),(x_{3},z_{2}),(x_{4},x_{0})\}\subseteq\Gamma_{2,s}$. Then $(x_{0},x_{1},x_{2},x_{3},x_{4},z_{2},x_{0})$ is an $s$-chain. By Lemma \ref{jcs2_3neq(3,3)} again, we get $\partial(x_{0},x_{3})=\partial(x_{3},x_{0})=2$ and $s=2$.

Note that there exists a vertex $w_{1}\in P_{(1,2),(1,2)}(x_{0},z_{2})\setminus\{z_{1}\}$. By $\Gamma_{1,2}\Gamma_{2,1}=\{\Gamma_{0,0},\Gamma_{3,3}\}$ and Lemma \ref{jcs2_3neq(3,3)}, we have $\Gamma_{1,2}(x_{0})=\{x_{1},z_{1},w_{1}\}$. Since $\Gamma_{1,2}\notin(\Gamma_{1,2})^{2}$, we get $z_{1}$ or $w_{1}\in P_{(1,2),(1,2)}(x_{0},x_{3})$. Observe $z_{3}\in P_{(2,1),(1,2)}(z_{1},x_{4})$ and $z_{2}\in P_{(1,2),(2,1)}(w_{1},x_{4})$, contrary to $\Gamma_{2,1}\Gamma_{1,2}=\Gamma_{1,2}\Gamma_{2,1}=\{\Gamma_{0,0},\Gamma_{3,3}\}$. Thus, $k_{3}=1$.$\qed$

\begin{lemma}\label{jcs2gammai2}
Let $(x_{0},x_{1},\ldots,x_{i})$ be an $s$-chain with $4\leq i\leq n-1$, where $n$ is the length of an $s$-line. Then $\Gamma_{i}=\Gamma_{\wz{\partial}(x_{0},x_{i})}$ and $k_{i}\in\{1,3\}$.
\end{lemma}
\textbf{Proof.}~Use induction on $i$ for $4\leq i \leq n-1$. Suppose that the lemma holds for $s$-chain $(x_{0},x_{1},\ldots,x_{j})$ with $0\leq j\leq i-1$. Write $\wz{\partial}(x_{0},x_{h})=\wz{a}_{h}$ for $0\leq h\leq i$.

\textbf{Case 1.}~$p_{(1,2),(3,3)}^{\wz{a}_{4}}\neq0$.

By Lemma \ref{jcs2gamma4}, $k_{3}=1$. Let $i=3l+r$, where $1\leq r\leq3$.

For any $(y_{0},y_{i})\in\Gamma_{\wz{a}_{i}}$, from $x_{r}\in P_{\wz{a}_{r},\wz{a}_{3l}}(x_{0},x_{i})$ and the inductive hypothesis, there exists an $s$-chain $(y_{0},y_{1},\ldots,y_{r})$ with $(y_{r},y_{i})\in\Gamma_{3l}$. Pick a vertex $y_{r+1}$ such that $y_{r+1}\in P_{(2,s),(2,1)}(y_{r-1},y_{r})$. By $\Gamma_{3l}=\Gamma_{\wz{a}_{3l}}\in(\Gamma_{\wz{a}_{3}})^{l}$ and Lemma \ref{jiben} (i), we get $k_{3l}=1$ and there exists an $s$-chain $(y_{r},y_{r+1},\ldots,y_{i})$. Then $(y_{0},y_{i})\in\Gamma_{i}$ and $\Gamma_{\wz{a}_{i}}\subseteq\Gamma_{i}$.

Conversely, let $(y_{0},y_{1},\ldots,y_{i})$ be an $s$-chain. By Lemma \ref{jiben} (iii), we get $\Gamma_{\wz{a}_{3l}}\Gamma_{\wz{a}_{r}}=\{\Gamma_{\wz{a}_{i}}\}$, which implies $\wz{\partial}(y_{0},y_{i})=\wz{a}_{i}$ from the inductive hypothesis. Hence, $\Gamma_{i}\subseteq\Gamma_{\wz{a}_{i}}$, so $\Gamma_{i}=\Gamma_{\wz{a}_{i}}$. By Lemma \ref{jiben} (i), if $r=1$ or $2$, then $k_{i}=3$; if $r=3$, then $k_{i}=1$.

\textbf{Case 2.} $p_{(1,2),(3,3)}^{\wz{a}_{4}}=0$.

Pick any $(y_{0},y_{i})\in \Gamma_{\wz{a}_{i}}$. Suppose $(y_{0},y_{i})\notin\Gamma_{i}$. By the inductive hypothesis and $x_{i-l}\in P_{\wz{a}_{i-l},\wz{a}_{l}}(x_{0},x_{i})$ with $l=1,2,\ldots,i-1$, there exist two $s$-chains $(y_{l,0}=y_{0},y_{l,1},\ldots,y_{l,i-l})$ and $(y_{l,i-l},y_{l,i-l+1},\ldots,y_{l,i}=y_{i})$ such that $\wz{\partial}(y_{l,i-l-1},y_{l,i-l+1})=(2,1)$. Since $\Gamma_{2,1}\Gamma_{1,2}=\Gamma_{1,2}\Gamma_{2,1}=\{\Gamma_{0,0},\Gamma_{3,3}\}$, $\{(y_{1,i-3},y_{i}),(y_{2,i-3},y_{i}),(y_{3,i-4},y_{3,i-1})\}\subseteq\Gamma_{3,3}$. By $\wz{\partial}(y_{1,i-2},y_{i})=(2,1)$, one has $y_{1,i-2}\neq y_{2,i-2}$ and $y_{1,1}\neq y_{2,1}$. From Lemma \ref{jcs2_3neq(3,3)}, we get $y_{3,i-3}\notin\{y_{1,i-3},y_{2,i-3}\}$ and $\Gamma_{1,2}(y_{0})=\{y_{1,1},y_{2,1},y_{3,1}\}$. Since $p_{(1,2),(3,3)}^{\wz{a}_{4}}=0$, we obtain $i\geq5$. If $k_{i-4}=1$, then $y_{1,i-4}=y_{4,i-4}$ and $y_{1,i-3}\in P_{(1,2),(3,3)}(y_{4,i-4},y_{i})$, a contradiction; if $k_{i-4}=3$, then $y_{4,i-4}\in\{y_{1,i-4},y_{2,i-4},y_{3,i-4}\}$, and $y_{3,i-1}\in P_{(3,3),(1,2)}(y_{4,i-4},y_{i})$ or $y_{h,i-3}\in P_{(1,2),(3,3)}(y_{4,i-4},y_{i})$ for some $h\in\{1,2\}$, a contradiction. Thus, $(y_{0},y_{i})\in\Gamma_{i}$ and $\Gamma_{\wz{a}_{i}}\subseteq\Gamma_{i}$.

Conversely, let $(y_{0},y_{1},\ldots,y_{i})$ be an $s$-chain. Suppose $\wz{\partial}(y_{0},y_{i})\neq\wz{\partial}(x_{0},x_{i})$. Since $\Gamma_{\wz{a}_{i}}\subseteq\Gamma_{i}$ and $x_{i}\in P_{\wz{a}_{i},\wz{a}_{l}^{*}}(x_{0},x_{i-l})$ with $l=1,2,3$, there exists a vertex $z_{l,i}\in P_{\wz{a}_{i},\wz{a}_{l}^{*}}(y_{0},y_{i-l})\cap\Gamma_{i}(y_{0})$ and $z_{l,i}\neq y_{i}$. Assume that $z\in P_{(1,2),(1,2)}(y_{i-2},z_{2,i})$. In view of $p_{(2,s),(2,1)}^{(1,2)}=1$ from Lemma \ref{jiben} (ii), we have $\wz{\partial}(y_{i-2},z_{1,i})=\wz{\partial}(y_{i-3},z)=(2,1)$ and $z_{1,i}\neq z_{2,i}$. By $\Gamma_{2,1}\Gamma_{1,2}=\Gamma_{1,2}\Gamma_{2,1}=\{\Gamma_{0,0},\Gamma_{3,3}\}$, we obtain $\wz{\partial}(y_{i-3},z_{l,i})=(3,3)$ for any $l=1,2$. Since $z_{3,i}\in \Gamma_{3}(y_{i-3})$, from Lemma \ref{jcs2_3neq(3,3)}, one gets $z_{3,i}\notin\{y_{i},z_{1,i},z_{2,i}\}$ and $\{y_{i},z_{1,i},z_{2,i},z_{3,i}\}\subseteq\Gamma_{i}(y_{0})$, contrary to $|\Gamma_{i}(y_{0})|\leq3$. Hence, $\wz{\partial}(y_{0},y_{i})=\wz{\partial}(x_{0},x_{i})$ and $\Gamma_{i}\subseteq\Gamma_{\wz{a}_{i}}$. Thus, $\Gamma_{i}=\Gamma_{\wz{a}_{i}}$.

At last, we will prove $k_{i}\in \{1,3\}$ for this case. Suppose $k_{i}=2$. For any vertex $x_{0}$, let $\Gamma_{1,2}(x_{0})=\{x_{1},y_{1},z_{1}\}$ and $(x_{0},x_{1},\ldots,x_{i})$, $(x_{0},y_{1},y_{2},\ldots,y_{i})$, $(x_{0},z_{1},z_{2},\ldots,z_{i})$ be three $s$-chains. Without loss of generality, we may assume $x_{i}=y_{i}$ and $z_{i}\neq x_{i}$. By Lemma \ref{jiben} (ii), we have $k_{\wz{a}_{i}}p_{\wz{a}_{i-1},(1,2)}^{\wz{a}_{i}}=k_{1,2}p_{\wz{a}_{i},\wz{a}_{i-1}^{*}}^{(1,2)}$,
which implies $p_{\wz{a}_{i-1},(1,2)}^{\wz{a}_{i}}=k_{\wz{a}_{i-1}}=3$. Then $\Gamma_{2,1}(x_{i})=\Gamma_{2,1}(z_{i})=\{x_{i-1},y_{i-1},z_{i-1}\}$. Since $p_{(2,s),(2,1)}^{(1,2)}=1$, we get $\{(z_{i-2},x_{i}),(y_{i-2},z_{i}),(x_{i-2},z_{i})\}\subseteq\Gamma_{2,1}$ and $k_{i-2}=3$. But $P_{\wz{a}_{i-2},(2,s)}(x_{0},x_{i})=\{x_{i-2},y_{i-2}\}$ and $P_{\wz{a}_{i-2},(2,s)}(x_{0},z_{i})=\{z_{i-2}\}$, a contradiction. Hence, $k_{i}\in\{1,3\}$.$\qed$

\begin{lemma}\label{jcs2coro}
Let $n$ be the length of an $s$-line.\vspace{-0.3cm}
\begin{itemize}
\item[${\rm(i)}$] If $(x_{0},x_{1},\ldots,x_{n})$ is an $s$-chain, then $x_{0}=x_{n}$.\vspace{-0.3cm}

\item[${\rm(ii)}$] For any $s$-chain $(y_{0},y_{1},\ldots,y_{l})$, $y_{i}=y_{j}$ if and only if $i\equiv j~({\rm mod}~n)$ with $0\leq i,j\leq l$.
\end{itemize}
\end{lemma}
\textbf{Proof.}~(i) By Lemmas \ref{jcs2_3neq(3,3)} and \ref{jcs2gammai2}, one has $\Gamma_{n-1}=\Gamma_{2,1}$ and $\Gamma_{n-2}=\Gamma_{s,2}$. Lemma \ref{jiben} (ii) implies $p_{(2,s),(2,1)}^{(1,2)}=1$. Since $x_{0},x_{n}\in P_{(2,s),(2,1)}(x_{n-2},x_{n-1})$, we get $x_{n}=x_{0}$.

(ii) is obvious by (i).$\qed$

\subsubsection{Construction of the digraph}

\begin{lemma}\label{jcs2line}
Let $L=(y_{0},y_{1},\ldots,y_{j})$ be an $s$-chain. Then the following hold:\vspace{-0.3cm}
\begin{itemize}
\item[${\rm(i)}$] For each $y_{0}'\in \Gamma_{1,2}(y_{0})\setminus\{y_{1}\}$, there exists a unique $s$-chain $L':=(y_{0}',y_{1}',\ldots,y_{j}')$
such that $(L, L')\in\wz{\Gamma}_{1,2}$.\vspace{-0.3cm}

\item[${\rm(ii)}$] There exists a unique $s$-chain $L''$ such that
$(L',L'')\in\wz{\Gamma}_{1,2}$ and $(L, L'')\in\wz{\Gamma}_{2,s}$.
\end{itemize}
\end{lemma}
\textbf{Proof.}~We define inductively $y_{i}'$ to be the unique element in $P_{(1,2),(2,1)}(y_{i-1}',y_{i})$ for $1\leq i\leq j$, which is well-defined by $A_{1,2}A_{2,1}=3I+A_{3,3}$ and $y_{i-1}'\neq y_{i}$ since $p_{(1,2),(1,2)}^{(2,s)}=1$. It remains to show that $L':=(y_{0}',y_{1}',\ldots,y_{j}')$ is an $s$-chain. Suppose for the contrary that $\wz{\partial}(y_{i}',y_{i+2}')=(2,1)$ for some $i\in\{0,1,\ldots,j-2\}$. By $\wz{\partial}(y_{i},y_{i+2})=(2,s)$, one gets $y_{i}\neq y_{i+2}'$. Since $\wz{\partial}(y_{i+1},y_{i+2}')=(2,1)$ from $p_{(1,2),(1,2)}^{(2,s)}=1$, we have $y_{i}', y_{i+1}\in P_{(1,2),(2,1)}(y_{i},y_{i+2}')$, contrary to $A_{1,2}A_{2,1}=3I+A_{3,3}$. This proves (i).

In view of $y_{0}'\neq y_{1}$ and $p_{(1,2),(1,2)}^{(2,s)}=1$, we get $\wz{\partial}(y_{0},y_{1}')=(2,1)$. Since $p_{(1,2),(s,2)}^{(2,1)}=1$ from Lemma \ref{jiben} (ii), there exists a unique vertex $y_{0}''\in P_{(1,2),(s,2)}(y_{0}',y_{0})\subseteq\Gamma_{1,2}(y_{0}')\setminus\{y_{1}'\}$ and a unique $s$-chain $L'' = (y_{0}'', y_{1}'',\ldots,y_{j}'')$ such that $(L',L'')\in\wz{\Gamma}_{1,2}$ by (i). Since $(y_{0},y_{0}',y_{0}'')$ is an $s$-chain, from (i) and the inductive hypothesis, we obtain $(y_{i},y_{i}'')\in\Gamma_{2,s}$ for $1\leq i\leq j$. Thus, (ii) is valid.$\qed$

By Lemma \ref{jcs2line}, there exists an $s$-plane $(L_{0},L_{1},\ldots,L_{n})$. For the remainder of this section, we assume that $L_{j}=(y(j,0),y(j,1),\ldots,y(j,n))$ and $y(i,j+1)\neq y(i+1,j)$ for $0\leq i,j\leq n$, where the coordinates could be read modulo $n$.

\begin{lemma}\label{jcs2k_3,3}
For $0\leq i,j\leq n-1$, the following hold:\vspace{-0.3cm}
\begin{itemize}
\item[${\rm(i)}$] $\partial(y(i,j),y(i-1,j-1))=1$ and $\wz{\partial}(y(i,j),y(i-2,j-2))=(2,s)$.\vspace{-0.3cm}

\item[${\rm(ii)}$] $\{y(i+1,j-1),y(i-1,j+1),y(i+2,j+1),y(i+1,j+2),y(i-2,j-1),y(i-1,j-2)\}\subseteq\Gamma_{3,3}(y(i,j))$.
\end{itemize}
\end{lemma}
\textbf{Proof.}~(i)~Since $y(i,j-1),y(i-1,j)\in P_{(1,2),(1,2)}(y(i-1,j-1),y(i,j))$ and $y(i,j-1)\neq y(i-1,j)$, we have $\partial(y(i,j),y(i-1,j-1))=1$ by $p_{(1,2),(1,2)}^{(2,s)}=1$.

Suppose $\wz{\partial}(y(i,j),y(i-2,j-2))=(2,1)$ for some
$i,j$. Since $y(i-1,j-2),y(i-2,j-1),y(i,j)\in P_{(1,2),(1,2)}(y(i-2,j-2),y(i-1,j-1))$ and $y(i-1,j-2)\neq y(i-2,j-1)$, one gets $y(i,j)\in\{y(i-1,j-2),y(i-2,j-1)\}$ by $p_{(1,2),(1,2)}^{(2,1)}=2$. In view of the symmetry, we may assume
$y(i,j)=y(i-1,j-2)$. It follows that $(2,s)=\wz{\partial}(y(i,j-2),y(i,j))=\wz{\partial}(y(i,j-2),y(i-1,j-2))=(2,1)$, a contradiction.

(ii)~Since $k_{1,2}=3$ and $\Gamma_{1,2}\Gamma_{2,1}=\{\Gamma_{0,0},\Gamma_{3,3}\}$, we have $\Gamma_{2,1}(y(h,l))\setminus\{y(i,j)\}\subseteq\Gamma_{3,3}(y(i,j))$ for any $(h,l)\in\{(i+1,j),(i,j+1),(i-1,j-1)\}$. Hence, (ii) is valid.$\qed$

For any vertex $y(i,j)$, there are three distinct $s$-lines beginning from $y(i,j)$. They are $(y(i,j),y(i+1,y),\ldots,y(i,j))$, $(y(i,j),y(i,j+1),\ldots,y(i,j))$ and $(y(i,j),y(i-1,j-1),\ldots,y(i,j))$. We consider the cardinality of $\{y(i,j)\mid 0\leq i,j\leq n-1\}$.

\begin{prop}\label{jcs2_n^2}
If $|\{y(i,j)\mid 0\leq i,j\leq n-1\}|=n^{2}$, then $\Gamma$ is isomorphic to one of the digraphs in Theorem {\rm\ref{Main} (vii)}.
\end{prop}
\textbf{Proof.}~Since $s\neq1$, we obtain $n\geq4$. If $3\mid n$, then $\wz{\partial}(y(0,0),y(1+n/3,0))=\wz{\partial}(y(0,0),y(1,n/3))=(1+n/3,-1+2n/3)$, contrary to  $y(1,n/3)\notin\Gamma_{1+n/3}(y(0,0))$. Thus, the desired result holds.$\qed$

Finally, we consider the case that $|\{y(i,j)\mid 0\leq i,j\leq n-1\}|\neq n^{2}$.

\begin{lemma}\label{jcs2kl}
For $1\leq l\leq n-1$, $k_{l}=1$ if and only if $3l\equiv0~({\rm mod}~n)$ and $y(a,b)=y(a+l,b+2l)=y(a+2l,b+l)$ for any $a,b\in\{0,1,\ldots,n-1\}$.
\end{lemma}
\textbf{Proof.}~"$\Longrightarrow$"~Since $k_{l}=1$, we have $y(a+l,b)=y(a,b+l)$, which implies that $y(a+l,b)=y(a-2l,b)$ by $|\Gamma_{l}(y(a-l,b+l))|=1$. Hence, $3l\equiv0$ (mod $n$). In view of $|\Gamma_{l}(y(a,b+2l))|=1$, one gets $y(a,b)=y(a+l,b+2l)=y(a+2l,b+l)$.

"$\Longleftarrow$"~Since $|\Gamma_{l}(y(a,b-l))|=1$, we obtain $k_{l}=1$.$\qed$

\begin{lemma}\label{jcs2_ki=kj}
If $y(a,b)=y(a+i,b+j)$ for some $a,b\in\{0,1,\ldots,n-1\}$ and $i,j\in\{1,2,\ldots,n-1\}$, then $k_{i}=k_{j}$.
\end{lemma}
\textbf{Proof.}~By Lemmas \ref{jcs2_3neq(3,3)} and \ref{jcs2gammai2}, we only need to prove that $k_{i}=1$ if and only if $k_{j}=1$. For any $l\in\{i,j\}$, from $|\Gamma_{l}(y(a,b-l))|\in\{1,3\}$ and Lemma \ref{jcs2kl}, $k_{l}=1$ if and only if $y(a,b)=y(a+l,b+2l)=y(a+2l,b+l)=y(a+i,b+j)$ and $3l\equiv0$ (mod $n$), if and only if $j\equiv2i~({\rm mod}~n)$ and $i\equiv2j~({\rm mod}~n)$. This completes the proof of the lemma.$\qed$

If $\Gamma_{l}=\Gamma_{\wz{f}}$ for some $\wz{f}\in\{\wz{i},\wz{j},\wz{h}\}$, then we may replace $\wz{f}$ with $l$ in the intersection number $p_{\wz{i},\wz{j}}^{\wz{h}}$ and the set $P_{\wz{i},\wz{j}}(x,y)$ for $x,y\in V\Gamma$.

\begin{lemma}\label{jcs2kiis1}
If $y(0,0)=y(i,j)$ with $1\leq i,j\leq n-1$ and $i\neq j$, then $k_{i}=1$.
\end{lemma}
\textbf{Proof.}~Suppose for the contrary that $k_{i}=3$.

\begin{step}\label{jcs2y=y1 1}
~{\rm Show that $n\nmid(h+l)$ and $y(0,0)=y(h-l,h)=y(2h,2l)$ if $k_{h}=3$ and $y(0,0)=y(h,l)$ with $n\nmid h$ and $n\nmid(h-l)$.}
\end{step}

By Lemma \ref{jcs2coro}, we may assume that $1\leq h,l\leq n-1$. Since $\Gamma_{h}(y(0,l))=\{y(h,l),y(0,h+l),y(-h,l-h)\}$, one gets $y(0,0)=y(h,l)\neq y(0,h+l)$, which implies $n\nmid(h+l)$. By $y(h,l)\in P_{0,n-l}(y(0,0),y(h,0))$, we get $\emptyset\neq P_{0,n-l}(y(h,l),y(h,h+l))\subseteq\{y(h+l,h+l),y(h,h+2l),y(h-l,h)\}$ and $\emptyset\neq P_{0,n-l}(y(h,l),y(2h,l))\subseteq\{y(2h,2l),y(2h+l,l),y(2h-l,0)\}$. Since $n\nmid(h+l)$ and $y(h,l)=y(0,0)$, one has $y(0,0)=y(h-l,h)$ and $y(0,0)\in\{y(2h,2l),y(2h-l,0)\}$. By $\Gamma_{h}(y(h-l,0))=\{y(2h-l,0),y(h-l,h),y(-l,-h)\}$, we obtain $y(0,0)=y(2h,2l)$.

\begin{step}\label{jcs2y=y1 3}
~{\rm Show that $\Gamma_{i+j}\neq\Gamma_{2n-i-j}$.}
\end{step}

Suppose for the contrary that $\Gamma_{i+j}=\Gamma_{2n-i-j}$. Observe that $p_{0,2n-i-j}^{i+j}=p_{0,i+j}^{i+j}\neq0$ and $P_{0,2n-i-j}(y(0,0),y(i+j,0))\subseteq\{y(i+j,i+j),y(2i+2j,0),y(0,-i-j)\}$. By Step \ref{jcs2y=y1 1}, we have $n\nmid(i+j)$ and $y(2i,2j)=y(i-j,i)$. Then $y(2i+2j,0)=y(0,0)$ and $i+j\equiv n/2~({\rm mod}~n)$. In view of Lemma \ref{jcs2_ki=kj}, we get $k_{n+i-j}=k_{i}=k_{j}=3$, which implies $y(2i,2j)=y(i-j,i)=y(i-j-i,i-j)=y(-2j,2i-2j)=y(2i,-4j)$. Hence, $6j\equiv0~({\rm mod}~n)$. Since $i+j\equiv n/2~({\rm mod}~n)$ and $j\neq i$, we obtain $(i,j)\in\{(n/3,n/6),(n/6,n/3),(5n/6,2n/3),(2n/3,5n/6)\}$. Then $y(i,j)=y(0,0)=y(2i,2j)=y(n/3,2n/3)$ or $y(2n/3,n/3)$. Since $i$ or $j\in\{n/3,2n/3\}$, by Lemma \ref{jcs2_ki=kj} and Step \ref{jcs2y=y1 1}, one has $y(0,0)=y(i,j)=y(4i,4j)=y(2n/3,n/3)$ or $y(n/3,2n/3)$, a contradiction.

\begin{step}\label{jcs2y=y1'}
~{\rm Show that $p_{(3,3),(3,3)}^{3}=p_{(3,3),(3,3)}^{n-3}=2$.}
\end{step}

Suppose $k_{3}=1$. By Lemma \ref{jcs2kl}, we have $n=9$ and $y(0,0)=y(3,6)=y(6,3)$, which imply $i\in\{1,2,4,5,7,8\}$ since $k_{i}=3$. Then $y(0,0)=y(1,l)$ for some $l$ by repeating using $k_{2i}=3$ and $y(0,0)=y(2i,2j)$ from Step \ref{jcs2y=y1 1}. In view of Lemma \ref{jcs2k_3,3} (ii), one gets $y(1,2),y(1,8)\in\Gamma_{3,3}(y(0,0))$ and $y(1,5)\in\Gamma_{3,3}(y(3,6))$. Since $\wz{\partial}(y(0,0),y(1,1))=(2,1)$, we obtain $l\in\{3,4,6,7\}$. By Step  \ref{jcs2y=y1 1}, one has $y(1-l,1)=y(1,l)=y(0,0)=y(6,3)=y(3,6)$, contrary to $n=9$. Hence, $k_{3}=3$.

Suppose $\Gamma_{3}=\Gamma_{n-3}$. By Lemma \ref{jcs2_3neq(3,3)}, we have $\Gamma_{3}=\Gamma_{2,2}$ and $s=2$. Pick a circuit $(x,y,z,w)$ such that $\wz{\partial}(x,z)=(2,2)$. Since $\Gamma_{1,2}\notin(\Gamma_{1,2})^{2}$, one gets $\wz{\partial}(y,w)=(2,2)$. Note that $(x,y,z,w,x)$ is an $2$-line, contrary to $\Gamma_{3}=\Gamma_{2,2}$. Thus, $\Gamma_{3}\neq\Gamma_{n-3}$

By Lemma \ref{jiben} (ii) and (iv), we get $k_{3}p_{(3,3),(3,3)}^{3}=k_{n-3}p_{(3,3),(3,3)}^{n-3}=k_{3,3}p_{(3,3),n-3}^{(3,3)}=k_{3,3}p_{(3,3),3}^{(3,3)}$ and $k_{3,3}\geq p_{(3,3),(3,3)}^{(3,3)}+p_{(3,3),3}^{(3,3)}+p_{(3,3),n-3}^{(3,3)}+p_{(3,3),(0,0)}^{(3,3)}$. Since $y(i-j+1,2-j),y(i-j+1,-j-1)\in P_{(3,3),(3,3)}(y(i-j+2,1-j),y(i-j,-j))$ and $y(0,0)\in P_{(3,3),(3,3)}(y(i-2,j-1),y(i+1,j-1))$ from Lemma \ref{jcs2k_3,3} (ii), one has $p_{(3,3),(3,3)}^{(3,3)}\geq2$ and $p_{(3,3),(3,3)}^{3}=p_{(3,3),(3,3)}^{n-3}\neq0$, which imply $p_{(3,3),3}^{(3,3)}=1$ and $p_{(3,3),(3,3)}^{3}=p_{(3,3),(3,3)}^{n-3}=2$.

\vspace{3ex}

Based on the above discussion, we reach a contradiction as follows.

For any $(h,l)\in\{(i+1,j-1),(i+1,j+2),(i-2,j-1)\}$, in view of Lemma \ref{jcs2k_3,3} (ii), we have $y(h,l)\in\Gamma_{3,3}(y(0,0))=\Gamma_{3,3}(y(i,j))$ and $y(h,0)\in P_{n+h,l}(y(0,0),y(h,l))$. Since $k_{3,3}=6$ and $k_{3,3}p_{n+h,l}^{(3,3)}=k_{n+h}p_{(3,3),2n-l}^{n+h}$ from Lemma \ref{jiben} (ii), we get $p_{(3,3),2n-l}^{n+h}=2$ by Lemmas \ref{jcs2_3neq(3,3)} and \ref{jcs2gammai2}. Hence, $|\Gamma_{3,3}(y(0,0))\cap\Gamma_{l}(y(h,0))|=2$. Since $\Gamma_{l}(y(h,0))=\{y(h,l),y(h+l,0),y(h-l,-l)\}$, we obtain $\{y(h+l,0),y(h-l,-l)\}\cap\Gamma_{3,3}(y(0,0))\neq\emptyset$.

By Step \ref{jcs2y=y1 3}, we have $y(i+j,0)\notin\Gamma_{3,3}(y(0,0))$ and $y(i-j+2,1-j)\in\Gamma_{3,3}(y(0,0))$. Suppose $y(i-j-1,a-j)\in\Gamma_{3,3}(y(0,0))$ for some $a\in\{1,-2\}$. From Lemma \ref{jcs2k_3,3} (ii), one gets $y(0,0),y(i-j,-j),y(i-j+1,1-j+a)\in P_{(3,3),(3,3)}(y(i-j+2,1-j),y(i-j-1,a-j))$. By Step \ref{jcs2y=y1 1}, we have $y(0,0)=y(i-j,i)$ and $y(0,0)\in\Gamma_{i+j}(y(i-j,-j))$, which imply $y(0,0)\neq y(i-j,-j)$. In view of $y(i-j+1,1-j+a)\in\Gamma_{3,3}(y(i-j,-j))$ and Step \ref{jcs2y=y1'}, we obtain $y(0,0)=y(i-j+1,1-j+a)$, which implies $\wz{\partial}(y(i-j,i),y(i-j,-j))=\wz{\partial}(y(0,0),y(i-j,-j))=\wz{\partial}(y(i-j+1,1-j+a),y(i-j,-j))=(3,3)$, contrary to Step \ref{jcs2y=y1 3}. Hence, $\{y(i+j+3,0),y(i+j-3,0)\}\subseteq\Gamma_{3,3}(y(0,0))$.

For any $(h',l')\in\{(i+1,j+2),(i-2,j-1)\}$, since $p_{0,2n-h'-l'}^{h'+l'}=p_{0,(3,3)}^{(3,3)}\neq0$, we have $y(0,0)\in\Gamma_{h'+l'}(y(h'+l',0))$. By $n\nmid(h'+l')$, $y(0,0)=y(2h'+2l',0)$. Therefore, $y(0,0)=y(2i+2j+6,0)=y(2i+2j-6,0)$. It follows that $12\equiv0~({\rm mod}~n)$. Since $n\geq4$, $n\in\{4,6,12\}$. By $n\nmid(i+j)$ from Step \ref{jcs2y=y1 1}, we get $i+j\equiv3~({\rm mod}~n)$ or $i+j\equiv n-3~({\rm mod}~n)$, contrary to $\{y(i+j+3,0),y(i+j-3,0)\}\subseteq\Gamma_{3,3}(y(0,0))$. Thus, $k_{i}=1$.$\qed$

\begin{prop}\label{jcs2_neq_n^2}
If $|\{y(i,j)\mid 0\leq i,j\leq n-1\}|\neq n^{2}$, then $3\mid n$, $n\geq9$ and $\Gamma\simeq{\rm Cay}(\mathbb{Z}_{n/3}\times\mathbb{Z}_{n},\{(0,1),(1,1),(-1,-2)\})$.
\end{prop}
\textbf{Proof.}~For some $a,b\in\{0,1,\ldots,n-1\}$, suppose $y(a,b)=y(a+i,b+j)$ with $i,j\in\{1,2\ldots,n-1\}$ and $i\neq j$. Since $y(a+i,b+j)\in P_{0,n-j}(y(a,b),y(a+i,b))$ and $y(a+i,b+j)\in P_{0,n+i-j}(y(a,b),y(a+i,b+i))$, one gets $y(0,0)\in P_{0,n-j}(y(0,0),y(i,0))\subseteq\{y(i,j),y(i+j,0),y(i-j,-j)\}$ and $y(0,0)\in P_{0,n+i-j}(y(0,0),y(-i,0))\subseteq\{y(j-2i,0),y(-i,j-i),y(-j,i-j)\}$. If $y(0,0)\in\{y(i,j),y(i-j,-j),y(-i,j-i),y(-j,i-j)\}$, from Lemmas \ref{jcs2_ki=kj} and \ref{jcs2kiis1}, then $k_{i}=k_{n-i}=1$ or $k_{j}=k_{n-j}=1$; if $y(0,0)=y(i+j,0)=y(j-2i,0)$, then $3i\equiv0~({\rm mod}~n)$ and $j\equiv2i~({\rm mod}~n)$, which imply that $y(a,b)=y(a+i,b+j)=y(a+i,b+2i)$. By Lemma \ref{jcs2_ki=kj} or $|\Gamma_{i}(a+2i,b)|\in\{1,3\}$, we have $k_{i}=1$. In view of Lemma \ref{jcs2kl}, one obtains $3\mid n$, $k_{n/3}=1$ and $(i,j)\in\{(n/3,2n/3),(2n/3,n/3)\}$.

Since $k_{1}=k_{2}=3$, $n\geq9$. By $k_{n/3}=1$, one gets $|\Gamma_{n/3}(y(c+2n/3,d))|=1$ for $0\leq c,d\leq n-1$, which implies $y(c,d)=y(c+n/3,d+2n/3)=y(c+2n/3,d+n/3)$ and $|\{y(i,j)\mid 0\leq i,j\leq n-1\}|=n^{2}/3$. Let $\sigma$ be the mapping from $\Gamma$ to $\textrm{Cay}(\mathbb{Z}_{n/3}\times\mathbb{Z}_{n},\{(0,1),(1,1),(-1,-2)\})$ such that $\sigma(a,b)=(a,a+b)$ for $0\leq a,b\leq n-1$. Routinely, $\sigma$ is a desired isomorphism. Thus, the desired result holds.$\qed$

\section{$p_{(1,g-1),(1,g-1)}^{(2,g-2)}=1$}

In this section, we prove Theorem \ref{Main} under the assumption that $p_{(1,g-1),(1,g-1)}^{(2,g-2)}=1$. We begin with the following lemma.

\begin{lemma}\label{jcs1AA}
The one of the following holds:\vspace{-0.3cm}
\begin{itemize}
\item [{\rm D1)}] $(A_{1,g-1})^{2}=A_{2,g-2}+A_{1,g-1}+3A_{2,s}$ and $k_{2,s}=1$.\vspace{-0.3cm}

\item [{\rm D2)}] $(A_{1,g-1})^{2}=A_{2,g-2}+A_{2,l}+3A_{2,s}$ and $k_{2,s}=1$.\vspace{-0.3cm}

\item [{\rm D3)}] $(A_{1,g-1})^{2}=A_{2,g-2}+3A_{2,s}$ and $k_{2,s}=1$.\vspace{-0.3cm}

\item [{\rm D4)}] $(A_{1,g-1})^{2}=A_{2,g-2}+2A_{2,s}$ and $k_{2,s}=3$.\vspace{-0.3cm}
\end{itemize}
Moreover, $p_{(1,g-1),(g-1,1)}^{\wz{h}}\leq1$ for $\wz{h}\neq(0,0)$.
\end{lemma}
\textbf{Proof.}~Observe $p_{(g-1,1),(1,g-1)}^{(0,0)}=3$ and $|\Gamma_{g-1,1}\Gamma_{1,g-1}|\geq2$. By Lemma \ref{jiben} (v), we have
$\sum_{\wz{i}\in\wz{\partial}(\Gamma)}p_{(1,g-1),(1,g-1)}^{\wz{i}}p_{(g-1,1),\wz{i}}^{(1,g-1)}=\sum_{\wz{j}\in\wz{\partial}(\Gamma)}p_{(g-1,1),(1,g-1)}^{\wz{j}}p_{\wz{j},(1,g-1)}^{(1,g-1)}>3$. It follows from Lemma \ref{jiben} (iv) that $p_{(1,g-1),(1,g-1)}^{\wz{i}}\geq2$ for some $\wz{i}\in\wz{\partial}(\Gamma)$.

\textbf{Case 1.} $\Gamma_{1,g-1}\in(\Gamma_{1,g-1})^{2}$.

By Lemma \ref{jiben} (i), one has $A_{1,g-1}A_{g-1,1}=3A_{0,0}+A_{1,g-1}+A_{g-1,1}$. In view of $p_{(1,g-1),(0,0)}^{(1,g-1)}=1$ and Lemma \ref{jiben} (iv), we obtain $p_{(1,g-1),(1,g-1)}^{(1,g-1)}=1$. Since $p_{(1,g-1),(1,g-1)}^{\wz{i}}\geq2$ for some $\wz{i}\in\wz{\partial}(\Gamma)$, from Lemma \ref{jiben} (vi), one gets $(A_{1,g-1})^{2}=A_{2,g-2}+A_{1,g-1}+3A_{2,s}$ and $k_{2,s}=1$. Thus, D1 holds.

\textbf{Case 2.} $\Gamma_{1,g-1}\notin(\Gamma_{1,g-1})^{2}$.

Since $p_{(1,g-1),(1,g-1)}^{\wz{i}}\geq2$ for some $\wz{i}\in\wz{\partial}(\Gamma)$, from Lemma \ref{jiben} (iii), we have $(\Gamma_{1,g-1})^{2}=\{\Gamma_{2,g-2},\Gamma_{2,s}\}$ or $(\Gamma_{1,g-1})^{2}=\{\Gamma_{2,g-2},\Gamma_{2,l},\Gamma_{2,s}\}$ with $g-2<l<s$.

We claim that $k_{2,s}\neq2$ and $k_{2,l}\neq1$. Suppose $k_{2,s}=2$. Let $(x,y_{0},z_{0},w)$ be a path such that $\wz{\partial}(x,z_{0})=(2,g-2)$ and $\wz{\partial}(y_{0},w)=(2,s)$. By Lemma \ref{jiben} (vi) and (ii), we obtain $p_{(1,g-1),(1,g-1)}^{(2,s)}=3$ and $p_{(g-1,1),(2,s)}^{(1,g-1)}=2$, which imply that there exist two vertices $z_{1}\in P_{(1,g-1),(1,g-1)}(y_{0},w)\setminus\{z_{0}\}$ and $y_{1}\in P_{(g-1,1),(2,s)}(z_{1},w)\setminus\{y_{0}\}$. Hence, $\wz{\partial}(x,z_{1})=(2,s)$ and $\partial(x,y_{1})=\partial(y_{1},z_{0})=1$, contrary to $p_{(1,g-1),(1,g-1)}^{(2,g-2)}=1$. Suppose $k_{2,l}=1$. Let $(x_{0},x_{1},x_{2},x_{3})$ be a path such that $\wz{\partial}(x_{0},x_{2})=(2,l)$ and $\wz{\partial}(x_{1},x_{3})=(2,s)$. It follows from Lemma \ref{jiben} (iii) that $|\Gamma_{2,l}\Gamma_{1,g-1}|=1$ and $\partial(x_{3},x_{0})=l-1$. Hence, $s=\partial(x_{3},x_{1})\leq\partial(x_{3},x_{0})+1=l$, contrary to $l<s$. Thus, our claim is valid.

If $|(\Gamma_{1,g-1})^{2}|=3$, from Lemma \ref{jiben} (i),(vi) and the claim, then $p_{(1,g-1),(1,g-1)}^{(2,l)}=1$, which implies that $k_{2,s}=1$ and $p_{(1,g-1),(1,g-1)}^{(2,s)}=3$; if $|(\Gamma_{1,g-1})^{2}|=2$, then $p_{(1,g-1),(1,g-1)}^{(2,s)}=3$ and $k_{2,s}=1$, or $p_{(1,g-1),(1,g-1)}^{(2,s)}=2$ and $k_{2,s}=3$. Thus, D2, D3 or D4 holds.

At last, we will prove $p_{(1,g-1),(g-1,1)}^{\wz{h}}\leq1$ with $\wz{h}\neq(0,0)$ for this case. We conclude $\sum_{\wz{j}\in\wz{\partial}(\Gamma)}p_{(g-1,1),(1,g-1)}^{\wz{j}}p_{\wz{j},(1,g-1)}^{(1,g-1)}=\sum_{\wz{i}\in\wz{\partial}(\Gamma)}p_{(1,g-1),(1,g-1)}^{\wz{i}}p_{(g-1,1),\wz{i}}^{(1,g-1)}=5$ from Lemma \ref{jiben} (i) and (ii). By Lemma \ref{jiben} (iii), one has $|\Gamma_{1,g-1}\Gamma_{g-1,1}|=2$ or $3$. If $\Gamma_{1,g-1}\Gamma_{g-1,1}=\{\Gamma_{0,0},\Gamma_{\wz{h}}\}$, from the commutativity of $\Gamma$, then $p_{\wz{h},(1,g-1)}^{(1,g-1)}k_{1,g-1}=p_{(g-1,1),(1,g-1)}^{\wz{h}}k_{\wz{h}}=6$, which implies that $p_{\wz{h},(1,g-1)}^{(1,g-1)}=2$ and $p_{(g-1,1),(1,g-1)}^{\wz{h}}=1$. If $\Gamma_{1,g-1}\Gamma_{g-1,1}=\{\Gamma_{0,0},\Gamma_{\wz{h}_{1}},\Gamma_{\wz{h}_{2}}\}$,
by Lemma \ref{jiben} (vi), then $p_{\wz{h}_{j},(1,g-1)}^{(1,g-1)}k_{1,g-1}=p_{(g-1,1),(1,g-1)}^{\wz{h}_{j}}k_{\wz{h}_{j}}=3$
for $j=1,2$, which implies $p_{\wz{h}_{j},(1,g-1)}^{(1,g-1)}=1$ and $p_{(g-1,1),(1,g-1)}^{\wz{h}_{j}}=1$.$\qed$

In the remaining of this section, we divide the proof into three subsections according to separate assumptions based on Lemma \ref{jcs1AA}.

\subsection{The case D1}

\begin{prop}\label{jcs1_Q_8}
In this case, $\Gamma$ is isomorphic to the digraph in Theorem {\rm\ref{Main} (ii)}.
\end{prop}
\textbf{Proof.}~By Lemma \ref{jiben} (i), we get $A_{1,g-1}A_{g-1,1}=3A_{0,0}+A_{1,g-1}+A_{g-1,1}$. Pick distinct vertices $x,y,z,w$ such that $\partial(x,y)=\partial(y,z)=\partial(x,z)=\partial(y,w)=1$ and $\wz{\partial}(x,w)=(2,g-2)$.
Observe $y\in P_{(g-1,1),(1,g-1)}(z,w)$. In view of the commutativity of $\Gamma$ and $p_{(1,g-1),(1,g-1)}^{(2,g-2)}=1$, one has $\partial(w,z)=1$. Since $z\in P_{(1,g-1),(g-1,1)}(x,w)$ and $\partial(x,w)=2$, we obtain $\partial(w,x)=1$ and $g=3$.

Choose a vertex $w'\in P_{(2,s),(2,1)}(x,y)$. Since $p_{(1,2),(1,2)}^{(1,2)}\neq0$, we obtain $w'\in P_{(1,2),(1,2)}(y,w)$, which implies $s=2$. By $k_{2,2}=1$ and Lemma \ref{jiben} (iii), one gets $\Gamma_{2,2}\Gamma_{1,2}=\{\Gamma_{2,1}\}$. Since $\Gamma_{1,2}\Gamma_{2,1}=\Gamma_{1,2}\Gamma_{2,1}$, we have $(\Gamma_{1,2})^{i}=\{\Gamma_{1,2},\Gamma_{2,1},\Gamma_{0,0},\Gamma_{2,2}\}$ for $i\geq3$. Then $\wz{\partial}(\Gamma)=\{(0,0),(1,2),(2,1),(2,2)\}$ and $|V\Gamma|=8$. By \cite{H}, the desired result follows.$\qed$

\subsection{The cases D2 and D3}

\begin{lemma}\label{pre-main}
If {\rm D2} or {\rm D3} holds, then $\Gamma_{2,g-2}\in\Gamma_{1,g-1}\Gamma_{g-1,1}$.
\end{lemma}
\textbf{Proof.}~The proof is rather long, and we shall prove it in Section 4.

\begin{prop}\label{jcs1pappus}
If {\rm D2} or {\rm D3} holds, then $\Gamma$ is isomorphic to the digraph in Theorem {\rm\ref{Main} (iv)}.
\end{prop}
\textbf{Proof.}~Pick a circuit $(x_{0,0},x_{1,0},\ldots,x_{g-1,0})$, where the first subscription of $x$ could be read modulo $g$. Since $p_{(2,s),(g-1,1)}^{(1,g-1)}=1$ from Lemma \ref{jiben} (ii), there exists a unique vertex $x_{i,1}\in P_{(2,s),(g-1,1)}(x_{i-1,0},x_{i,0})$. The fact that $p_{(1,g-1),(1,g-1)}^{(2,s)}=3$ implies that $(x_{0,1},x_{1,1},\ldots,x_{g-1,1})$ is a minimal circuit.

Since $p_{(1,g-1),(g-1,1)}^{(2,g-2)}\neq0$ from Lemma \ref{pre-main}, there exists a vertex $x_{0}$ such that $x_{0}\in P_{(1,g-1),(g-1,1)}(x_{0,0},x_{2,0})$. The fact  $p_{(1,g-1),(1,g-1)}^{(2,s)}=3$ and $(x_{0,0},x_{1,1}),(x_{2,0},x_{3,1})\in\Gamma_{2,s}$ imply that $\Gamma_{1,g-1}(x_{0,0})=\{x_{1,0},x_{0,1},x_{0}\}$ and $\Gamma_{1,g-1}(x_{2,0})=\{x_{2,1},x_{3,0},x_{0}\}$. Similarly, there exists a vertex $x_{1}\in P_{(1,g-1),(g-1,1)}(x_{0,1},x_{2,1})$ such that $(x_{1,1},x_{1})\notin\Gamma_{2,s}$. By $\Gamma_{1,g-1}\notin(\Gamma_{1,g-1})^{2}$, we obtain $\partial(x_{2,0},x_{0,0})\neq1$ and $g\geq4$.

Suppose that $(A_{1,g-1})^{2}=A_{2,g-2}+A_{2,l}+3A_{2,s}$ with $g-2<l<s$. By $p_{(2,l),(g-1,1)}^{(1,g-1)}\neq0$, we obtain $(x_{1,0},x_{0})\in\Gamma_{2,l}$. Since $x_{0,0}\in P_{(g-1,1),(1,g-1)}(x_{1,0},x_{0})$, from the commutativity of $\Gamma$, one has $\Gamma_{2,l},\Gamma_{l,2},\Gamma_{2,g-2},\Gamma_{g-2,2}\in\Gamma_{1,g-1}\Gamma_{g-1,1}$. In view of Lemma \ref{jiben} (iii), we get $(2,l)\in\{(g-2,2),(l,2)\}$, contrary to $g-2<l=2$. Then $(A_{1,g-1})^{2}=A_{2,g-2}+3A_{2,s}$. By Lemma \ref{jiben} (i), we have $k_{2,g-2}=6$ and $\Gamma_{1,g-1}\Gamma_{g-1,1}=\{\Gamma_{0,0},\Gamma_{2,g-2}\}$. Thus, $g=4$.

Pick a vertex $x_{2}$ such that $\Gamma_{1,3}(x_{2,1})=\{x_{3,1},x_{1},x_{2}\}$. Since $x_{1},x_{3,1}\in\Gamma_{2,2}(x_{1,1})$ and $p_{(2,s),(3,1)}^{(1,3)}\neq0$, we have $(x_{1,1},x_{2})\in\Gamma_{2,s}$ and $(x_{2,0},x_{2})\in\Gamma_{2,2}$ from $k_{2,s}=1$. The fact that $p_{(1,3),(3,1)}^{(2,2)}\neq0$ implies $x_{0}$ or $x_{3,0}\in P_{(1,3),(3,1)}(x_{2,0},x_{2})$. By $p_{(1,3),(1,3)}^{(2,s)}=3$, one gets $x_{0}\in P_{(1,3),(1,3)}(x_{0,0},x_{1,1})$. From $\partial(x_{2},x_{1,1})=s>2$, $(x_{2},x_{0})$ is not an arc and $\partial(x_{2},x_{3,0})=1$. Since $(x_{2},x_{3,0},x_{0,0},x_{0,1},x_{1,1})$ is a path, we have $s\leq4$. If $s=3$, by $k_{2,3}=1$ and Lemma \ref{jiben} (iii), then $(\Gamma_{2,3})^{2}=\{\Gamma_{3,1}\}$ or $\Gamma_{2,3}\Gamma_{2,2}=\{\Gamma_{3,1}\}$, contrary to Lemma \ref{jiben} (i). Thus, $s=4$.

Since $(\Gamma_{1,3})^{2}=\{\Gamma_{2,2},\Gamma_{2,4}\}$, we have $\partial(x_{0,0},x_{2,1})=3$. By $k_{2,4}=1$ and Lemma \ref{jiben} (iii), one gets $\partial(x_{2,1},x_{0,0})=3$ and $\Gamma_{2,4}\Gamma_{1,3}=\{\Gamma_{3,3}\}$, which imply $k_{3,3}=3$ from Lemma \ref{jiben} (i). Since $\partial(x_{3,0},x_{0,0})=\partial(x_{0,0},x_{0})=1$, one gets $\Gamma_{2,2}\Gamma_{1,3}=\{\Gamma_{3,1},\Gamma_{1,3},\Gamma_{3,3}\}$ and $(\Gamma_{1,3})^{3}=\{\Gamma_{1,3},\Gamma_{3,1},\Gamma_{3,3}\}$.

Note that $x_{1,1}\in P_{(2,4),(2,4)}(x_{0,0},x_{2})$ and $\partial(x_{2},x_{0,0})=2$. By $k_{2,4}=1$ and Lemma \ref{jiben} (i), $\wz{\partial}(x_{0,0},x_{2})=(4,2)$. The fact $x_{1}\neq x_{1,1}$ implies $(x_{0,0},x_{1})\in\Gamma_{2,2}$. Since $x_{0,1}\in P_{(1,3),(3,1)}(x_{3,1},x_{0,0})$, one gets $\wz{\partial}(x_{0,0},x_{3,1})=(2,2)$. By $\Gamma_{1,3}(x_{2,1})=\{x_{1},x_{2},x_{2,2}\}$, one has $\Gamma_{3,3}\Gamma_{1,3}=\{\Gamma_{2,2},\Gamma_{4,2}\}$ and $(\Gamma_{1,3})^{4}=\{\Gamma_{2,2},\Gamma_{4,2},\Gamma_{0,0},\Gamma_{2,4}\}$ from $\Gamma_{3,1}\Gamma_{1,3}=\Gamma_{1,3}\Gamma_{3,1}$. By Lemma \ref{jiben} (iii), one obtains $\Gamma_{4,2}\Gamma_{1,3}=\{\Gamma_{3,1}\}$. Then $(\Gamma_{1,3})^{2n-1}=\{\Gamma_{1,3},\Gamma_{3,1},\Gamma_{3,3}\}$ and $(\Gamma_{1,3})^{2n}=\{\Gamma_{2,2},\Gamma_{4,2},\Gamma_{0,0},\Gamma_{2,4}\}$ for $n\geq2$. Therefore,  $\wz{\partial}(\Gamma)=\{(0,0),(1,3),(2,2),(2,4),(3,1),(3,3),(4,2)\}$ and $|V\Gamma|=18$. By \cite{H}, the desired result holds.$\qed$

\subsection{The case D4}

In this case, by Lemma \ref{jiben} (vi), we have $k_{2,g-2}=3$, which implies $p_{(2,g-2),(g-1,1)}^{(1,g-1)}=1$ and $p_{(2,s),(g-1,1)}^{(1,g-1)}=2$ from Lemma \ref{jiben} (ii).

Note that any $(g-2)$-line is a minimal circuit and the length of a $(g-2)$-line is $g$. As an application of $(g-2)$-line, we give a construction of the digraph $\Gamma$.

\subsubsection{Construction of $(g-2)$-planes}

\begin{lemma}\label{jcs1_1}
Let $x$ be a vertex and $\Gamma_{1,g-1}(x)=\{y_{0},y_{1},y_{2}\}$.\vspace{-0.3cm}
\begin{itemize}
\item[${\rm(i)}$]~If there exist distinct vertices $z$ and $w$ such that $z\in P_{(2,s),(g-1,1)}(x,y_{0})$ and $w\in P_{(1,g-1),(g-1,1)}(y_{0},y_{2})$~(See Figure {\rm1}), then $\partial(y_{1},z)=1$.\vspace{-0.3cm}

\item[${\rm(ii)}$]~If there exist three vertices $z_{0},z_{1},z_{2}$ such that $\Gamma_{g-1,1}(x)=\{z_{0},z_{1},z_{2}\}$ and $\wz{\partial}(z_{0},y_{0})=(2,g-2)$, then $\{(z_{0},y_{1}),(z_{0},y_{2}),(z_{1},y_{0}),(z_{2},y_{0})\}\subseteq\Gamma_{2,s}$.
\end{itemize}
\end{lemma}

\begin{figure}[hbt]
  \centering
   \scalebox{0.75}{\includegraphics{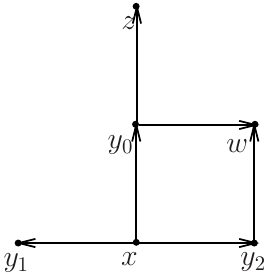}}
   \caption{ Lemma \ref{jcs1_1} (i).}
\end{figure}
\textbf{Proof.}~(i)~By Lemma \ref{jcs1AA}, we get $\partial(y_{2},z)\neq1$. Since $p_{(1,g-1),(1,g-1)}^{(2,s)}=2$, one has $y_{1}\in P_{(1,g-1),(1,g-1)}(x,z)$

(ii) Since $p_{(g-1,1),(2,s)}^{(1,g-1)}=p_{(2,s),(g-1,1)}^{(1,g-1)}=2$, one has $z_{1},z_{2}\in P_{(g-1,1),(2,s)}(x,y_{0})$ and $y_{1},y_{2}\in P_{(2,s),(g-1,1)}(z_{0},x)$.$\qed$

\begin{lemma}\label{jcs1coro}
Let $L=(y_{0},y_{1},\ldots,y_{g-1},y_{g})$ be a $(g-2)$-line. Then the following hold:\vspace{-0.3cm}
\begin{itemize}
\item[${\rm(i)}$] $\hat{L}:=(y_{1},y_{2},\ldots,y_{g},y_{1})$ is a $(g-2)$-line and $(L,\hat{L})\in\wz{\Gamma}_{1,g-1}$.\vspace{-0.3cm}

\item[${\rm(ii)}$] For each $y_{0}'\in \Gamma_{1,g-1}(y_{0})\setminus\{y_{1}\}$, there exists a unique $(g-2)$-line $L':=(y_{0}',y_{1}',\ldots,y_{g}')$
such that $(L, L')\in\wz{\Gamma}_{1,g-1}$.\vspace{-0.3cm}

\item[${\rm(iii)}$] There exists a unique $(g-2)$-line $L_{1}$ such that
$(L',L_{1})\in\wz{\Gamma}_{1,g-1}$ and $(L, L_{1})\in\wz{\Gamma}_{2,g-2}$.
\end{itemize}
\end{lemma}
\textbf{Proof.}~(i) is obvious.

(ii) Let $\Gamma_{1,g-1}(y_{0})=\{y_{1},y_{0}',y_{0}''\}$. Observe $y_{0}\in P_{(g-1,1),(1,g-1)}(y_{0}',y_{1})$ and $y_{0}\in P_{(g-1,1),(1,g-1)}(y_{0}'',y_{1})$. By Lemma \ref{jcs1AA} and the commutativity of $\Gamma$, there only exist two distinct vertices $y_{1}'\in P_{(1,g-1),(g-1,1)}(y_{0}',y_{1})$ and $y_{1}''\in P_{(1,g-1),(g-1,1)}(y_{0}'',y_{1})$ from $(A_{1,g-1})^{2}=A_{2,g-2}+2A_{2,s}$. Since $p_{(2,g-2),(g-1,1)}^{(1,g-1)}=1$, we may assume that $L':=(y_{0}',y_{1}',\ldots,y_{g-1}',y_{0}')$ is the unique $(g-2)$-line containing the arc $(y_{0}',y_{1}')$ and $(y_{0}'',y_{1}'',\ldots,y_{g-1}'',y_{0}'')$ is the unique $(g-2)$-line containing the arc $(y_{0}'',y_{1}'')$. It suffices to show that $(L,L')\in\wz{\Gamma}_{1,g-1}$. By induction, we only need to prove $\partial(y_{2},y_{2}')=1$.

Pick vertices $z_{0},z_{1},z_{1}',z_{2}$ such that $\Gamma_{1,g-1}(y_{0}'')=\{z_{0},z_{2},y_{1}''\}$ and $\partial(y_{1}',z_{1})=\partial(y_{1}'',z_{1}')=1$ with $z_{1}\neq y_{2}'$ and $z_{1}'\neq y_{2}''$. Since $p_{(2,g-2),(g-1,1)}^{(1,g-1)}=1$ and $p_{(2,s),(g-1,1)}^{(1,g-1)}=2$, we may assume
$$\tilde{\partial}(y_{i},z_{i})=\tilde{\partial}(y_{1},z_{1}')=(2,s)\ \textrm{and}\ \tilde{\partial}(y_{0},z_{2})=(2,g-2),$$
where $i=0,1$. If $z_{1}=y_{2}''$ and $z_{1}'=y_{2}'$, then $y_{2}',y_{2}''\in P_{(1,g-1),(g-1,1)}(y_{1}',y_{1}'')$, contrary to Lemma \ref{jcs1AA}. Without loss of generality, we may assume $z_{1}\neq y_{2}''$.

Suppose $\partial(y_{2},y_{2}')\neq1$.  Since $p_{(1,g-1),(1,g-1)}^{(2,s)}=2$ and $\wz{\partial}(y_{1},y_{2}')=(2,s)$ from Lemma \ref{jcs1_1} (ii), we have $\partial(y_{1}'',y_{2}')=1$. Note that $y_{1}''\in P_{(1,g-1),(g-1,1)}(y_{0}'',y_{1})$ and $y_{2}'\in P_{(1,g-1),(g-1,1)}(y_{1}',y_{1}'')$. By Lemma \ref{jcs1_1} (i) and $\wz{\partial}(y_{i},z_{i})=(2,s)$ for $i=0,1$, one gets $\partial(y_{0}',z_{0})=\partial(y_{2},z_{1})=1$. Observe $y_{1}',y_{1}''\in P_{(2,s),(1,g-1)}(y_{0},y_{2}')$ and $y_{0}'\in P_{(1,g-1),(2,g-2)}(y_{0},y_{2}')$. In view of the commutativity of $\Gamma$, we obtain $(y_{0}'',y_{2}')\in\Gamma_{2,s}$ and $y_{2}'\neq y_{2}''$. Since $y_{2}'\in P_{(1,g-1),(g-1,1)}(y_{1}',y_{1}'')$ and $\wz{\partial}(y_{1},y_{2}'')=(2,s)$ from Lemma \ref{jcs1_1} (ii), one has $\partial(y_{2},y_{2}'')=1$ by Lemma \ref{jcs1_1} (i). In view of $y_{1}''\in P_{(1,g-1),(g-1,1)}(y_{0}'',y_{1})$ and Lemma \ref{jcs1AA}, we get $y_{1}'\neq z_{0}$. Since $p_{(1,g-1),(1,g-1)}^{(2,s)}=2$ and $p_{(1,g-1),(1,g-1)}^{(2,g-2)}=1$, one obtains $z_{2}\in P_{(1,g-1),(1,g-1)}(y_{0}'',y_{2}')$.

By Lemma \ref{jcs1_1} (ii), one has $y_{0}',y_{1}\in P_{(1,g-1),(2,s)}(y_{0},z_{1})$. In view of $y_{2}\in P_{(2,g-2),(1,g-1)}(y_{0},z_{1})$, we get $\tilde{\partial}(y_{0}'',z_{1})=(2,g-2)$ from the commutativity of $\Gamma$. Since $y_{1}'\in P_{(g-1,1),(1,g-1)}(y_{2}',z_{1})$ and $y_{2}\in P_{(g-1,1),(1,g-1)}(y_{2}'',z_{1})$, by Lemma \ref{jcs1AA}, one has $\partial(z_{2},z_{1})\neq1$ and $\partial(y_{1}'',z_{1})\neq1$. Thus, $\partial(z_{0},z_{1})=1$.

Since $z_{1}\neq y_{2}''$ and $z_{0},y_{1}''\in\Gamma_{2,s}(y_{0})$, one has $|\{z_{1},y_{2}'',y_{0}\}|=3$. If $g=3$, then $y_{2}'',z_{1},y_{0}\in P_{(1,g-1),(1,g-1)}(y_{2},y_{0}'')$, contrary to $(A_{1,g-1})^{2}=A_{2,g-2}+2A_{2,s}$. Suppose $g>3$. Since $y_{1}',z_{0}\in\Gamma_{2,s}(y_{0})$ and $y_{2}\in\Gamma_{2,g-2}(y_{0})$, we have $\Gamma_{g-1,1}(z_{1})=\{z_{0},y_{2},y_{1}'\}$. By $y_{0}''\in P_{(1,g-1),(2,g-2)}(y_{g-1,1}'',z_{1})$ and $\tilde{\partial}(y_{g-1}'',z_{0})=(2,s)$ from Lemma \ref{jcs1_1} (ii), we get $y_{1}'$ or $y_{2}\in\Gamma_{2,g-2}(y_{g-1}'')$. In view of $\tilde{\partial}(y_{2},y_{3}'')=(2,s)$ and $g>3$, one obtains $\partial(y_{g-1}'',y_{2})\neq2$, which implies $\tilde{\partial}(y_{g-1}'',y_{1}')=(2,g-2)$. Then there exists a vertex $w\in P_{(1,g-1),(1,g-1)}(y_{g-1}'',y_{1}')$. Observe $y_{0}\in P_{(g-1,1),(1,g-1)}(y_{0}'',y_{1})$ and $y_{0}\in P_{(g-1,1),(1,g-1)}(y_{0}',y_{0}'')$. By Lemma \ref{jcs1AA}, $w\notin\{y_{0}',y_{1}\}$. Since $p_{(2,g-2),(g-1,1)}^{(1,g-1)}=1$, we get $\tilde{\partial}(w,z_{1})=(2,g-2)$ and $\tilde{\partial}(y_{g-1}'',z_{1})=(3,g-3)$, contrary to $\tilde{\partial}(y_{g-1}'',z_{0})=(2,s)$.

Hence, $\partial(y_{2},y_{2}')=1$. This proves (ii).

(iii) Since $p_{(2,g-2),(g-1,1)}^{(1,g-1)}=1$, there exists a vertex $w_{0}\in P_{(2,g-2),(g-1,1)}(y_{0},y_{0}')$ and a unique $(g-2)$-line $L_{1}=(w_{0},w_{1},\ldots,w_{g-1},w_{0})$ such that $(L',L_{1})\in\wz{\Gamma}_{1,g-1}$ by (ii). It suffices to show that $(y_{i},w_{i})\in\Gamma_{2,g-2}$ for $1\leq i\leq g-1$. By induction, we only need to prove $\wz{\partial}(y_{1},w_{1})=(2,g-2)$.

Suppose $\wz{\partial}(y_{1},w_{1})=(2,s)$. Assume that $\Gamma_{1,g-1}(y_{0})=\{y_{1},y_{0}',y_{0}''\}$. By (ii), there exists a unique $(g-2)$-line $L''=(y_{0}'',y_{1}'',\ldots,y_{g-1}'',y_{0}'')$ such that $(L,L'')\in\wz{\Gamma}_{1,g-1}$. Since $p_{(1,g-1),(1,g-1)}^{(2,g-2)}=1$, we have $w_{0}\neq y_{1}'$ and $(y_{0},y_{1}'),(y_{0},y_{1}''),(y_{0}',w_{1})\in\Gamma_{2,s}$. By $p_{(1,g-1),(1,g-1)}^{(2,s)}=2$ and $y_{2}\in\Gamma_{2,g-2}(y_{0})$, one gets $\Gamma_{1,g-1}(y_{1})=\{y_{1}',y_{1}'',y_{2}\}$. In view of $y_{2}'\in\Gamma_{2,g-2}(y_{0}')$, $w_{1}\neq y_{2}'$. Since $y_{2}'\in P_{(1,g-1),(g-1,1)}(y_{1}',y_{2})$, from Lemma \ref{jcs1_1} (i), one gets $y_{1}''\in P_{(1,g-1),(1,g-1)}(y_{1},w_{1})$. By $y_{1}'',y_{2}\in P_{(1,g-1),(1,g-1)}(y_{1},y_{2}'')$, one has $w_{1}\neq y_{2}''$ and $\wz{\partial}(y_{0}'',w_{1})=(2,s)$ from Lemma \ref{jcs1_1} (ii). But $y_{0}'',y_{1},y_{0}'\in P_{(1,g-1),(2,s)}(y_{0},w_{1})$ and $w_{0}\in P_{(2,g-2),(1,g-1)}(y_{0},w_{1})$, contrary to the commutativity of $\Gamma$.$\qed$

By Lemma \ref{jcs1coro} (ii) and (iii), there exists a $(g-2)$-plane $(L_{0},L_{1},\ldots,L_{g})$.

\subsubsection{Construction of the digraph}

For each $\wz{h}\in\wz{\partial}(\Gamma)$, we define a relation $\overline{\Gamma}_{\wz{h}}$ on the set of chains of $(g-2)$-lines as follows. For any two chains of $(g-2)$-lines $P=(L_{0},L_{1},\ldots,L_{j})$ and $P'=(L_{0}',L_{1}',\ldots,L_{j'}')$, $(P, P')\in\overline{\Gamma}_{\wz{h}}$ if and only if $j=j'$ and $(L_{i},L_{i}')\in \wz{\Gamma}_{\wz{h}}$ for $0\leq i\leq j$.

\begin{lemma}\label{jcs1coro2}
Suppose that $P=(L_{0},L_{1},\ldots,L_{g-1},L_{g})$ is a $(g-2)$-plane with $L_{j}=(y_{j,0},y_{j,1},\ldots,y_{j,g})$ and $y_{h,l+1}\neq y_{h+1,l}$ for $0\leq j\leq g$ and $0\leq h,l\leq g-1$. The following hold:\vspace{-0.3cm}
\begin{itemize}
\item[${\rm(i)}$] Both $\tilde{P}:= (L_{1},L_{2},\ldots,L_{g},L_{1})$ and $\hat{P}:=(\hat{L}_{0},\hat{L}_{1},\ldots,
\hat{L}_{g})$ are $(g-2)$-planes such that $(P,\tilde{P})\in\overline{\Gamma}_{1,g-1}$ and $(P,\hat{P})\in\overline{\Gamma}_{1,g-1}$.\vspace{-0.3cm}

\item[${\rm(ii)}$] There exists a unique $(g-2)$-plane $P'=(L_{0}',L_{1}',\ldots,L_{g-1}',L_{g}')$ with $P'\notin\{\tilde{P},\hat{P}\}$ and $(P, P')\in\overline{\Gamma}_{1,g-1}$.\vspace{-0.3cm}

\item[${\rm(iii)}$] There exists a unique $(g-2)$-plane $P''$ with
$(P',P'')\in\overline{\Gamma}_{1,g-1}$ and $(P,P'')\in\overline{\Gamma}_{2,g-2}$.
\end{itemize}
\end{lemma}
\textbf{Proof.}~Without loss of generality, we may assume that the first and second subscriptions of $y$ could be read modulo $g$. Note that (i) is obvious.

(ii) Let $y_{0,0}'\in\Gamma_{1,g-1}(y_{0,0})\setminus\{y_{1,0},y_{0,1}\}$. By Lemma \ref{jcs1coro} (ii), there exists a unique $(g-2)$-line $L_{0}'=(y_{0,0}',y_{0,1}'\ldots,y_{0,g}')$  such that $(L_{0},L_{0}')\in\wz{\Gamma}_{1,g-1}$. It follows from Lemma \ref{jcs1_1} (ii) that $\wz{\partial}(y_{0,0},y_{0,1}')=(2,s)$. In view of $p_{(2,s),(g-1,1)}^{(1,g-1)}=2$ and Lemma \ref{jcs1coro} (ii), there exist a unique vertex $y_{1,0}'$ and a unique $(g-2)$-line $L_{1}'=(y_{1,0}',y_{1,1}',\ldots,y_{1,g}')$ such that $P_{(2,s),(g-1,1)}(y_{0,0},y_{0,0}')=\{y_{1,0}',y_{0,1}'\}$ and $(L_{0}',L_{1}')\in\wz{\Gamma}_{1,g-1}$. By Lemma \ref{jcs1coro} (iii), we can construct a unique $(g-2)$-plane $P'=(L_{0}',L_{1}',\ldots,L_{g-1}',L_{g}')$ .

Since $(y_{0,0},y_{1,0}')\in\Gamma_{2,s}$ and $y_{0,1}'\in P_{(1,g-1),(g-1,1)}(y_{0,1},y_{0,0}')$, from Lemma \ref{jcs1_1} (i), we get $\partial(y_{1,0},y_{1,0}')=1$. Since $\wz{\partial}(y_{1,g-1},y_{1,1})=\wz{\partial}(y_{0,g-1}',y_{0,1}')=(2,g-2)$, from Lemma \ref{jcs1_1} (ii), one has $y_{0,0}\notin\{y_{1,g-1},y_{0,g-1}'\}$. By $y_{0,0}'\in P_{(1,g-1),(g-1,1)}(y_{0,0},y_{0,g-1}')$ and Lemma \ref{jcs1AA}, we obtain $\Gamma_{1,g-1}(y_{0,g-1})=\{y_{0,0},y_{0,g-1}',y_{1,g-1}\}$. Since $\Gamma_{1,g-1}(y_{0,g-1})=\Gamma_{2,s}(y_{1,0}')$, from the commutativity of $\Gamma$, we have $\wz{\partial}(y_{0,g-1},y_{1,g-1}')=(2,s)$. In view of Lemma \ref{jcs1_1} (i) again, $\partial(y_{1,g-1},y_{1,g-1}')=1$. Similarly, $(L_{1},L_{1}')\in\wz{\Gamma}_{1,g-1}$. By induction, (ii) is valid.

(iii) Let $P'=(L_{0}',L_{1}',\ldots,L_{g-1}',L_{g}')$ with $L_{i}'=(y_{i,0}',y_{i,1}',\ldots,y_{i,g}')$ for $0\leq i\leq g$. By (ii), there exists a unique $(g-2)$-plane $P''$ such that $(P',P'')\in\overline{\Gamma}_{1,g-1}$. Assume that $P''=(L_{0}'',L_{1}'',\ldots,L_{g-1}'',L_{g}'')$ with $L_{i}''=(y_{i,0}'',y_{i,1}'',\ldots,y_{i,g}'')$ and $y_{0,0}''\notin\{y_{1,0}',y_{0,1}'\}$. Since $(A_{1,g-1})^{2}=A_{2,g-2}+2A_{2,s}$, from the inductive hypothesis, one has $y_{i,h}''\notin\{y_{i+1,h}',y_{i,h+1}'\}$ for $0\leq i,h\leq g$. In view of Lemma \ref{jcs1_1} (ii), we get $y_{i+1,h}',y_{i,h+1}'\in P_{(2,s),(g-1,1)}(y_{i,h},y_{i,h}')$, which implies $\wz{\partial}(y_{i,h},y_{i,h}'')=(2,g-2)$ by $p_{(2,s),(g-1,1)}^{(1,g-1)}=2$. This completes the proof of (iii).$\qed$

Let $P_{0},P_{1},\ldots,P_{j}$ be $(g-2)$-planes. We say that $(P_{0},P_{1},\ldots,P_{j})$ is a {\em chain of
$(g-2)$-planes} if $(P_{i},P_{i+1})\in\overline{\Gamma}_{1,g-1}$ for $0\leq i\leq j-1$, and $(P_{i'},P_{i'+2})\in\overline{\Gamma}_{2,g-2}$
for $0\leq i'\leq j-2$. In particular, we say that a chain of $(g-2)$-planes
$(P_{0},P_{1},\ldots,P_{j})$ is a {\em $(g-2)$-cube} if $P_{0}=P_{j}$.

In view of Lemma \ref{jcs1coro2} (ii) and (iii), there exists a $(g-2)$-cube $(P_{0},P_{1},\ldots,P_{g})$. In the remaining of this subsection, we assume that $P_{i}=(L_{i,0},L_{i,1},\ldots,L_{i,g})$ and $L_{i,j}=(y(i,j,0),y(i,j,1),\ldots,y(i,j,g))$ with $\Gamma_{1,g-1}(y(i,j,k))=\{y(i+1,j,k),y(i,j+1,k),y(i,j,k+1)\}$ for $0\leq i,j,k\leq g$, where the coordinates are taken modulo $g$.

\begin{lemma}\label{jcs1plane}
The following hold:\vspace{-0.3cm}
\begin{itemize}
\item[${\rm(i)}$] For any fixed $i$, $y(i,j,k)$s are pairwise distinct.\vspace{-0.3cm}

\item[${\rm(ii)}$] For any fixed $j$, $y(i,j,k)$s are pairwise distinct.\vspace{-0.3cm}

\item[${\rm(iii)}$] For any fixed $k$, $y(i,j,k)$s are pairwise distinct.
\end{itemize}
\end{lemma}
\textbf{Proof.}~We only need to prove (i). Suppose $y(i,0,0)=y(i,j,k)$ with $1\leq j,k\leq g-1$. Since $\partial(y(i,j,0),y(i,0,0))=\partial(y(i,j,0),y(i,j,k))=k$ and $(y(i,j,0),y(i,j+1,0),\ldots,y(i,0,0))$ is a path of length $g-j$, we obtain $g-j\geq k$. Observe that $C=(y(i,0,0),y(i,1,0),\ldots,y(i,j,0),y(i, j,1),\ldots,y(i,j,k-1))$ is a circuit of length $j+k$. Then $C$ is a minimal circuit and $\wz{\partial}(y(i,j-1,0),y(i,j,1))=(2,g-2)$, contrary to $\wz{\partial}(y(i,j-1,0),y(i,j,1))=(2,s)$ by Lemma \ref{jcs1_1} (ii).$\qed$

\begin{prop}\label{jcs1_z_3^3}
In this case, $\Gamma$ is isomorphic to the digraph in Theorem {\rm\ref{Main} (v)}.
\end{prop}
\textbf{Proof.}~We prove it step by step.

\begin{stepp}\label{jcs1spacelem 1}
~{\rm Show that $g<i+j+k<2g$ if $y(a,b,c)=y(a+i,b+j,c+k)$ for some $a,b,c\in\{0,1,\ldots,g-1\}$ and $i,j,k\in\{1,2,\ldots,g-1\}$.}
\end{stepp}

Note that $C=(y(a,b,c),y(a+1,b,c),\ldots,y(a+i,b,c),y(a+i,b+1,c),\ldots,y(a+i,b+j,c),y(a+i,b+j,c+1),\ldots,y(a+i,b+j,c+k-1))$ and $C'=(y(a+i,b+j,c+k),y(a+i+1,b+j,c+k),\ldots,y(a,b+j,c+k),y(a,b+j+1,c+k),\ldots,y(a,b,c+k),y(a,b,c+k+1),\ldots,y(a,b,c-1))$ are two circuits of length $i+j+k$ and $3g-i-j-k$, respectively. Hence, $g\leq i+j+k\leq2g$. By Lemma \ref{jcs1_1} (ii), one gets $(y(a+i-1,b,c),y(a+i,b+1,c)),(y(a-1,b+j,c+k),y(a,b+j+1,c+k))\in\Gamma_{2,s}$, which implies that neither $C$ nor $C'$ is a minimal circuit. Thus, $g<i+j+k<2g$.

\begin{stepp}\label{jcs1spacelem 2}
~{\rm Show that $y(d,e,f)\in\{y(d+m,e+l,f+h),y(d+m,e+h,f+l)\}$ with $\{h,l,m\}=\{i,j,k\}$ for any $y(d,e,f)$ if $y(a,b,c)=y(a+i,b+j,c+k)$ for some $a,b,c\in\{0,1,\ldots,g-1\}$ and $g\nmid i$.}
\end{stepp}

By Lemma \ref{jcs1plane}, we may assume $0<i,j,k<g$. We only need to consider the case that $m=k$. Since $\wz{\partial}(y(a+i,b+j,c),y(a,b,c))=\wz{\partial}(y(a+i,b+j,c),y(a+i,b+j,c+k))=(k,g-k)$ and $y(a,b+j,c)\in P_{(g-i,i),(g-j,j)}(y(a+i,b+j,c),y(a,b,c))$, we get $p_{(g-i,i),(g-j,j)}^{(k,g-k)}\neq0$, which implies $y(d,e,f)\in\Gamma_{i,g-i}(x)$ for some $x\in\Gamma_{j,g-j}(y(d+k,e,f))=\{y(d+k+j,e,f),y(d+k,e+j,f),y(d+k,e,f+j)\}$. By Lemma \ref{jcs1plane} and Step \ref{jcs1spacelem 1}, one obtains $y(d,e,f)\in\{y(d+k,e+j,f+i),y(d+k,e+i,f+j)\}$. The desired result holds.

\begin{stepp}\label{jcs1spacelem 22}
~{\rm Show that $i=j=k=g/2$ if $y(0,0,0)=y(i,j,k)$ and $|\{i,j,k\}|\neq3$ for some $i,j,k\in\{1,2,\ldots,g-1\}$.}
\end{stepp}

By Step \ref{jcs1spacelem 2}, $y(i,j,k)=y(0,0,0)\in\{y(j,k,i),y(j,i,k)\}\cap\{y(k,j,i),y(k,i,j)\}$. Since $|\{i,j,k\}|\neq3$, from Lemma \ref{jcs1plane}, one gets $i=j=k$. Suppose $i\neq\frac{g}{2}$. In view of Step \ref{jcs1spacelem 2}, one gets $y(0,0,0)=y(2i,2i,2i)=y(2i-g,2i-g,2i-g)$. By Step \ref{jcs1spacelem 1}, we have $i<g/3$ or $2g/3<i$, contrary to $g/3<i<2g/3$ from $y(0,0,0)=y(i,i,i)$. Hence, $i=g/2$.

\begin{stepp}\label{jcs1spacelem 222}
~{\rm Show that $y(i,j,k)=y(2i,2j,2k)=y(i+k,j+i,k+j)$ if $y(0,0,0)=y(i,j,k)$ with $g\nmid i$ and $g\nmid(i-j)$.}
\end{stepp}

By Step \ref{jcs1spacelem 22}, we have $g\nmid(j-k)$ and $g\nmid(i-k)$. It follows from Step \ref{jcs1spacelem 2} that $y(i,j,k)\in\{y(2i,j+k,k+j),y(2i,2j,2k)\}$ and $y(i,j,k)\in\{y(i+k,2j,i+k),y(i+k,j+i,k+j)\}$. If $y(i,j,k)=y(2i,j+k,k+j)=y(i+k,j+i,k+j)$ or $y(i,j,k)=y(2i,2j,2k)=y(i+k,2j,i+k)$, then $g\mid(k-i)$, a contradiction. Suppose $y(0,0,0)=y(i,j,k)=y(2i,j+k,k+j)=y(i+k,2j,i+k)$. By Step \ref{jcs1spacelem 22}, we obtain $i+k\equiv2j\equiv0$ or $g/2~({\rm mod}~g)$, and $j+k\equiv2i\equiv0$ or $g/2~({\rm mod}~g)$. Since $g\nmid(i-j)$, one has
$i+k\equiv2j\equiv0~({\rm mod}~g)$ and $j+k\equiv2i\equiv g/2~({\rm mod}~g)$, or $i+k\equiv2j\equiv g/2~({\rm mod}~g)$ and $j+k\equiv2i\equiv0~({\rm mod}~g)$. By $g\nmid i$ and Lemma \ref{jcs1plane}, we have $g\nmid j$, which implies $j\equiv g/2~({\rm mod}~g)$ or $i\equiv g/2~({\rm mod}~g)$. Then $g\mid k$, a contradiction. Thus, $y(i,j,k)=y(2i,2j,2k)=y(i+k,j+i,k+j)$.

\vspace{3ex}

Based on the above discussion, we complete the proof of this proposition as follows.

First, we will show that all vertices $y(d,e,f)$ with $0\leq d,e,f\leq g-1$ are distinct. Suppose $y(a,b,c)=y(a+i,b+j,c+k)$ for some $a,b,c\in\{0,1,\ldots,g-1\}$ and $i,j,k\in\{1,2,\ldots,g-1\}$. By Step \ref{jcs1spacelem 2}, $y(0,0,0)\in\{y(i,j,k),y(i,k,j)\}$. Without loss of generality, we may assume $y(0,0,0)=y(i,j,k)$.

Suppose $|\{i,j,k\}|=3$. By Lemma \ref{jcs1plane} and Step \ref{jcs1spacelem 222}, we have $y(0,0,0)=y(i+k,j+i,k+j)=y(2i+2k,2j+2i,2k+2j)$. It follows from Step \ref{jcs1spacelem 2} that $y(i+k,j+i,k+j)\in\{y(i+j+k,i+j+k,i+j+k),y(i+j+k,2i+j,2k+j)\}$ and $y(i+j+k,2i+j,2k+j)\in\{y(i+j+2k,2i+2j,2k+j+i),y(i+j+2k,3i+j,2k+2j)\}$. Suppose $y(2i+2k,2j+2i,2k+2j)=y(i+k,j+i,k+j)=y(i+j+k,2i+j,2k+j)$. By Lemma \ref{jcs1plane}, we obtain $i=j$, contrary to $|\{i,j,k\}|=3$. Hence, $y(0,0,0)=y(i+k,j+i,k+j)=y(i+j+k,i+j+k,i+j+k)=y(i+j+k-g,i+j+k-g,i+j+k-g)$. It follows from the Steps \ref{jcs1spacelem 1} and \ref{jcs1spacelem 22} that $i+j+k=3g/2$ and $y(0,0,0)=y(i+k,j+i,k+j)=y(g/2-j,g/2-k,g/2-i)$. Hence, $g/2-j+g/2-k+g/2-i\equiv0~({\rm mod}~g)$, contrary to $g\nmid(g/2-j+g/2-k+g/2-i)$ from Step \ref{jcs1spacelem 1}. Hence, $|\{i,j,k\}|\neq3$.

By Steps \ref{jcs1spacelem 2} and \ref{jcs1spacelem 22}, we get $y(d,e,f)=y(d+g/2,e+g/2,f+g/2)$ for $0\leq d,e,f\leq g-1$. Then $|\{y(d,e,f)\mid0\leq d,e,f\leq g-1\}|=g^{3}/2$. Observe $y(2,1,0),y(1,1,1)\in\Gamma_{3,3g/2-3}(y(0,0,0))$. But we have $y(2,0,0)\in P_{(2,g-2),(1,g-1)}(y(0,0,0),y(2,1,0))$ and $P_{(2,g-2),(1,g-1)}(y(0,0,0),y(1,1,1)))=\emptyset$,
a contradiction. Thus, all vertices $y(d,e,f)$ with $0\leq d,e,f\leq g-1$ are distinct.

Suppose $g>3$. Observe $y(2,2,0),y(3,1,0)\in\Gamma_{4,2g-4}(y(0,0,0))$. But $y(2,0,0)\in P_{(2,g-2),(2,g-2)}(y(0,0,0),y(2,2,0))$ and $P_{(2,g-2),(2,g-2)}(y(0,0,0),y(3,1,0))=\emptyset$,
a contradiction. The desired result holds.$\qed$

Combining Propositions \ref{jcs2h=3(1)}, \ref{jcs2h=3}, \ref{jcs2p=3}, \ref{jcs2_n^2}, \ref{jcs2_neq_n^2}, \ref{jcs1_Q_8}, \ref{jcs1pappus} and \ref{jcs1_z_3^3}, we complete the proof of Theorem \ref{Main}.

\section{Proof of Lemma \ref{pre-main}}

In this section, we will prove Lemma \ref{pre-main} by contradiction. Suppose $\Gamma_{2,g-2}\notin\Gamma_{1,g-1}\Gamma_{g-1,1}$. Since $p_{(1,g-1),(1,g-1)}^{(2,s)}=3$, from Lemma \ref{jiben} (ii), we get $p_{(2,s),(g-1,1)}^{(1,g-1)}=1$.

Let $t=\min\{i\mid p_{(i+1,g-i-1),(g-1,1)}^{(i,g-i)}=1\}$. Note that $1\leq t\leq g-2$. If $1\leq j<t$, by $p_{(2,s),(g-1,1)}^{(1,g-1)}=1$, then $p_{(j+1,g-j-1),(g-1,1)}^{(j,g-j)}=2$; if $j\geq t$, then $p_{(j+1,g-j-1),(g-1,1)}^{(j,g-j)}=1$.

Pick a minimal circuit $(x_{0,0},x_{1,0},\ldots,x_{g-1,0})$, where the first subscription of $x$ could be read modulo $g$. Since $k_{2,s}=1$, there exists a unique vertex $x_{i,j}$ such that $\wz{\partial}(x_{i-1,j-1},x_{i,j})=(2,s)$ for $0\leq i\leq g-1$ and $1\leq j\leq g+2$. By $p_{(1,g-1),(1,g-1)}^{(2,s)}=3$,  $(x_{0,j},x_{1,j},\ldots,x_{g-1,j})$ is a minimal circuit. Let $x_{i,j}'$ denote the vertex such that $\Gamma_{1,g-1}(x_{i,j})=\Gamma_{g-1,1}(x_{i+1,j+1})=\{x_{i+1,j},x_{i,j+1},x_{i,j}'\}$ for $0\leq i,j\leq g$.

By $p_{(2,s),(g-1,1)}^{(1,g-1)}=1$, we get $(x_{i,j}',x_{i+1,j+1}')\in\Gamma_{2,s}$ for $0\leq i,j\leq g-1$. The fact that $p_{(t+1,g-t-1),(g-1,1)}^{(t,g-t)}=1$ implies $(x_{i,j},x_{i+t,j}'),(x_{i,j}',x_{i+t+1,j+1})\notin\Gamma_{t+1,g-t-1}$. Since $p_{(1,g-1),(1,g-1)}^{(1,g-1)}=0$, we obtain $\partial(x_{0,0},x_{2,1})=3$. It follows from Lemma \ref{jiben} (iii) that $\partial(x_{2,1},x_{0,0})=s-1$ and $\Gamma_{2,s}\Gamma_{1,g-1}=\{\Gamma_{3,s-1}\}$. Thus, $\partial(x_{0,0},x_{2,2})=4$.

In the following, we give some useful results for distance of two vertices. Write $\wz{f}_{t}:=\wz{\partial}(x_{0,1},x_{t,1}')$ and $s_{i}:=\partial(x_{1+i,1+i},x_{0,0})$ for $0\leq i\leq3$. Note that $s_{0}=s$.

\begin{lemma}\label{jcs1_neq_2,s}
 Let $(x_{3,3}=y_{0,j},y_{1,j},\ldots,y_{s_{j},j}=x_{2-j,2-j})$ be a shortest path for $0\leq j\leq2$. If $\partial(x_{0,0},x_{3+j,1+j})=4+2j$ for some $j\in\{0,1,2\}$, then we have:\vspace{-0.3cm}
\begin{itemize}
\item[${\rm(i)}$] $s_{j+1}\neq s_{j}-2$ and $\wz{\partial}(y_{i,j},y_{i+2,j})\neq(2,s)$ for $0\leq i\leq s_{j}-2$.\vspace{-0.3cm}

\item[${\rm(ii)}$] $g\geq2t$ and $|\Gamma_{t,g-t}\Gamma_{1,g-1}|=3$.\vspace{-0.3cm}

\item[${\rm(iii)}$] $\wz{\partial}(y_{i-t,j},y_{i+1,j})=\wz{f}_{t}$ for some $i\in\{t,t+1,\ldots,s_{j}-1\}$.
\end{itemize}
\end{lemma}
\textbf{Proof.}~(i) Note that $\partial(x_{0,0},x_{2+j,1+j})=3+2j$ and $\wz{\partial}(x_{0,0},x_{1+j,1+j})=(2+2j,s_{j})$. Since $k_{2,s}=1$, from Lemma \ref{jiben} (i), we have $k_{2+2j,s_{j}}=1$. By Lemma \ref{jiben} (iii), we get $x_{2+j,1+j},x_{1+j,2+j},x_{1+j,1+j}'\in\Gamma_{3+2j,s_{j}-1}(x_{0,0})$. Hence, $\partial(x_{0,0},x_{2+j,2+j})=4+2j$.

Observe $\partial(x_{3+j,1+j},x_{0,0})\geq s_{j}-2$ and $s_{j+1}\geq s_{j}-2$. Suppose $s_{j+1}=s_{j}-2$. By Lemma \ref{jiben} (i), one gets $k_{4+2j,s_{j+1}}=1$, which implies $\partial(x_{3+j,2+j},x_{0,0})=s_{j}-3$ from Lemma \ref{jiben} (iii). Then $s_{j}-2\leq\partial(x_{3+j,1+j},x_{0,0})\leq1+\partial(x_{3+j,2+j},x_{0,0})\leq s_{j}-2$. But $x_{1+j,1+j}\in P_{(2+2j,s_{j}),(2,s)}(x_{0,0},x_{2+j,2+j})$ and $P_{(2+2j,s_{j}),(2,s)}(x_{0,0},x_{3+j,1+j})=\emptyset$, a contradiction. Therefore, $s_{j+1}\neq s_{j}-2$.

Suppose $\wz{\partial}(y_{i,j},y_{i+2,j})=(2,s)$ for some $i\in\{0,1,\ldots,s_{j}-2\}$. By the commutativity of $\Gamma$, we may assume $i=0$. Since $k_{2,s}=1$, one has $y_{2,j}=x_{4,4}$ and $\partial(x_{4,4},x_{2-j,2-j})=s_{j}-2$. It follows from Lemma \ref{jiben} (iii) that $\wz{\partial}(x_{2-j,2-j},x_{4,4})=\wz{\partial}(x_{0,0},x_{2+j,2+j})=(4+2j,s_{1+j})$, contrary to $s_{1+j}\neq s_{j}-2$. Thus, (i) is valid.

Let $\alpha_{j}=\max\{c\mid\wz{\partial}(y_{0,j},y_{c,j})=(c,g-c)\}$. Note that $\alpha_{j}<s_{j}$ and there exists a minimal circuit $(y_{0,j}=v_{g},y_{1,j},\ldots,y_{\alpha_{j},j}=v_{\alpha_{j}},v_{\alpha_{j}+1},v_{\alpha_{j}+2},\ldots,v_{g-1})$.

(ii) and (iii)\quad First, we prove $g\geq2t$ and consider the case $t\geq2$. Since $p_{(2,s),(g-1,1)}^{(1,g-1)}=1$, from (i), we get $\alpha_{j}\geq t$. Suppose $\alpha_{j}\geq g-t+1$. By $p_{(g-t,t),(1,g-1)}^{(g-t+1,t-1)}=p_{(t,g-t),(g-1,1)}^{(t-1,g-t+1)}=2$, there exists a vertex $v_{\alpha_{j}-1}\in P_{(g-t,t),(1,g-1)}(y_{\alpha_{j}-g+t-1,j},y_{\alpha_{j},j})\setminus\{y_{\alpha_{j}-1,j}\}$. Since $v_{\alpha_{j}-1},y_{\alpha_{j}-1,j}\in\Gamma_{g-2,2}(y_{\alpha_{j}+1,j})\cap\Gamma_{g-2,2}(v_{\alpha_{j}+1})$ and $p_{(2,s),(g-1,1)}^{(1,g-1)}=1$, we get $y_{\alpha_{j}+1,j}=v_{\alpha_{j}+1}$ and $\alpha_{j}\neq g-1$, contrary to the fact that $(y_{0,j},y_{1,j},\ldots,y_{\alpha_{j},j}=v_{\alpha_{j}},y_{\alpha_{j}+1,j}=v_{\alpha_{j}+1},v_{\alpha_{j}+2},\ldots,v_{g-1})$ is a minimal circuit. Thus, $\alpha_{j}\leq g-t$ and $g\geq2t$.

By (i), $\wz{\partial}(y_{\alpha_{j}-t,j},y_{\alpha_{j},j})=(t,g-t)$. The fact that $\partial(x_{0,1},x_{t-1,1})=t-1$ and $\wz{\partial}(x_{t-1,1},x_{t,2})=(2,s)$ imply $\wz{\partial}(y_{\alpha_{j}-t,j},y_{\alpha_{j}+1,j})\neq\wz{\partial}(x_{0,1},x_{t,2})$. Since  $y_{\alpha_{j}+1,j}\neq v_{\alpha_{j}+1}$ and $p_{(t+1,g-t-1),(g-1,1)}^{(t,g-t)}=1$, one gets $\wz{\partial}(y_{\alpha_{j}-t,j},y_{\alpha_{j}+1,j})\neq(t+1,g-t-1)$. Hence, $\wz{\partial}(y_{\alpha_{j}-t,j},y_{\alpha_{j}+1,j})=\wz{f}_{t}$. Thus, (ii) and (iii) hold.$\qed$

\begin{lemma}\label{jcs1d=4}
The following hold:\vspace{-0.3cm}
\begin{itemize}
\item[${\rm(i)}$] $g\geq5$ and $\partial(x_{0,0},x_{3,1})=4$.\vspace{-0.3cm}

\item[${\rm(ii)}$] If $t\geq2$, then $g\geq6$, $\partial(x_{0,0},x_{3,2})=\partial(x_{0,0},x_{4,1})=5$ and $\partial(x_{0,0}',x_{3,1}')+2=\partial(x_{0,0},x_{4,2})$.
\end{itemize}
\end{lemma}
\textbf{Proof.}~For $0\leq i,j\leq g$, we claim that $(x_{i,j}',x_{i+t+1,j}')\in\Gamma_{t+1,g-t-1}$ if $g\geq2t+2$. Since $x_{i+t+1,j}\in P_{(t+1,g-t-1),(1,g-1)}(x_{i,j},x_{i+t+1,j}')$, from the commutativity of $\Gamma$, we have $x_{i,j}'$ or $x_{i,j+1}\in P_{(1,g-1),(t+1,g-t-1)}(x_{i,j},x_{i+t+1,j}')$. The fact $p_{(t+1,g-t-1),(1,g-1)}^{(t+2,g-t-2)}=p_{(g-t-1,t+1),(g-1,1)}^{(g-t-2,t+2)}=1$ implies that $x_{i+t+1,j}'\notin P_{(t+1,g-t-1),(1,g-1)}(x_{i,j+1},x_{i+t+2,j+1})$. Thus, our claim is valid.

(i) Suppose $g=3$. Since $p_{(2,1),(2,1)}^{(1,2)}=p_{(1,2),(1,2)}^{(2,1)}=1$, $t=1$. By Lemma \ref{jcs1_neq_2,s} (ii), we have $(A_{1,2})^{2}=A_{2,1}+A_{2,l}+3A_{2,s}$ with $1<l<s\leq4$. From $x_{0,1}\in P_{(1,2),(2,1)}(x_{0,0},x_{2,1})$ and $p_{(1,2),(1,2)}^{(2,s)}=3$, we get $\wz{\partial}(x_{0,0},x_{2,1})=(3,s-1)=(3,3)$. Since $\Gamma_{2,4}\Gamma_{1,2}=\{\Gamma_{3,3}\}$, we obtain $x_{0,2}\in P_{(2,4),(1,2)}(x_{2,1},x_{0,0})$. By $(x_{0,0},x_{1,1})\in\Gamma_{2,4}$ and $p_{(2,l),(2,1)}^{(1,2)}\neq0$, one has $(x_{0,0},x_{0,1}')\in\Gamma_{2,l}$. The fact that $(x_{2,0}',x_{0,1}')\in\Gamma_{2,4}$ implies $(x_{2,1},x_{0,1}')\in\Gamma_{2,1}$, contrary to $(x_{2,1},x_{0,1}')\notin\Gamma_{2,1}$. Hence, $g\neq3$.

Suppose $\partial(x_{0,0},x_{3,1})<4$. Since $\partial(x_{0,0},x_{2,1})=\partial(x_{0,0},x_{3,0})=3$, $\partial(x_{0,0},x_{2,0}')<3$. By $p_{(1,g-1),(g-1,1)}^{(2,g-2)}=0$, one has $\partial(x_{0,0},x_{2,0}')=2$. The fact that $p_{(1,g-1),(g-1,1)}^{(1,g-1)}=0$ and $\partial(x_{0,1},x_{2,0}')\geq2$ imply $\partial(x_{0,0}',x_{2,0}')=1$. Similarly, $\partial(x_{2,0}',x_{4,0}')=1$. Then $g>4$. Since $P_{(3,g-3),(g-1,1)}(x_{0,0},x_{2,0})=\{x_{3,0}\}$, from the claim, we get $t=2$.

Since $p_{(1,g-1),(1,g-1)}^{(2,s)}=3$, we have $x_{2,0}'\in P_{(1,g-1),(1,g-1)}(x_{0,0}',x_{1,1}')$. The fact $\partial(x_{0,0}',x_{4,0})\neq1$ implies $(x_{0,0}',x_{4,0}')\notin\Gamma_{2,s}$ and $x_{1,1}'\neq x_{4,0}'$. By $p_{(1,g-1),(1,g-1)}^{(1,g-1)}=0$, one gets $\Gamma_{1,g-1}(x_{2,0}')=\{x_{1,1}',x_{4,0}',x_{3,1}\}$ and $x_{2,0}'\neq x_{1,1}$. Hence, $(x_{0,0},x_{2,0}')\in\Gamma_{2,g-2}$ and $\Gamma_{2,g-2}\in\Gamma_{2,g-2}\Gamma_{1,g-1}$. Observe $x_{1,1}\in P_{(2,s),(1,g-1)}(x_{0,0},x_{1,1}')$ and $x_{1,1}\in P_{(2,s),(2,g-2)}(x_{0,0},x_{3,1})$. By $x_{1,1}'\in P_{(1,g-1),(g-1,1)}(x_{1,1},x_{2,0}')$ and $p_{(1,g-1),(g-1,1)}^{(1,g-1)}=p_{(1,g-1),(g-1,1)}^{(2,g-2)}=0$, we obtain $x_{3,1},x_{1,1}'\notin\Gamma_{2,g-2}(x_{0,0})$, which implies $\wz{\partial}(x_{0,0},x_{4,0}')=(2,g-2)$. Since $\partial(x_{0,1},x_{4,0}')\geq4$, one has $\partial(x_{1,0},x_{4,0}')=1$. From $p_{(1,g-1),(1,g-1)}^{(2,s)}=3$, we get $x_{4,0}'\in P_{(1,g-1),(1,g-1)}(x_{1,0},x_{2,1})$. Note that $x_{2,1}\in P_{(1,g-1),(g-1,1)}(x_{2,0},x_{4,0}')$, contrary to $p_{(1,g-1),(g-1,1)}^{(2,g-2)}=0$ or $p_{(1,g-1),(1,g-1)}^{(2,s)}=3$. Then $\partial(x_{0,0},x_{3,1})=4$.

Suppose $g=4$. Note that $3\leq s\leq6$ and $\wz{\partial}(x_{0,0},x_{2,1})=(3,s-1)$. Since $p_{(1,3),(1,3)}^{(2,s)}=3$ and $p_{(1,3),(1,3)}^{(1,3)}=0$, we have $s\neq3$. By Lemma \ref{jcs1_neq_2,s} (i), $s_{1}>s-2$ and $(x_{3,2},x_{0,0})$ is not an arc, which imply $s\neq4$ from $\Gamma_{2,s}\Gamma_{1,3}=\{\Gamma_{3,s-1}\}$. Since $x_{1,1}\in P_{(2,s),(2,2)}(x_{0,0},x_{3,1})$ and $P_{(2,s),(2,2)}(x_{0,0},x_{2,2})=\emptyset$, one gets $\partial(x_{2,2},x_{0,0})\neq\partial(x_{3,1},x_{0,0})$. If $\partial(x_{3,1},x_{0,0})=4$, then $x_{0,2}\in P_{(2,s),(2,2)}(x_{3,1},x_{0,0})$, which implies $\partial(x_{2,2},x_{0,0})<4$, contrary to $s_{1}>s-2$; if $\partial(x_{3,1},x_{0,0})=3$, then $s=5$, which implies $x_{0,2}\in P_{(2,5),(1,3)}(x_{3,1},x_{0,0})$, contrary to $g=4$. Thus, (i) holds.

(ii) Observe $x_{1,1}\in\Gamma_{2,s}(x_{0,0})\cap\Gamma_{g-2,2}(x_{3,1})\cap\Gamma_{g-2,2}(x_{2,1}')$. It follows from (i) and Lemma \ref{jiben} (iii) that $\partial(x_{0,0},x_{3,2})=5$.

Suppose $g=5$. (i) and Lemma \ref{jcs1_neq_2,s} (ii) imply $t=2$.  By Lemma \ref{jcs1_neq_2,s} (i), $s-2<\partial(x_{2,2},x_{0,0})\leq6$. Hence, $4\leq s\leq7$ and $\Gamma_{4,s-2}\notin(\Gamma_{2,s})^{2}$. In view of Lemma \ref{jiben} (iii), one gets $\Gamma_{2,s}\Gamma_{2,3}=\{\Gamma_{4,s-2}\}$ and $(x_{0,0},x_{3,1})\in\Gamma_{4,s-2}$. Since $p_{(1,4),(1,4)}^{(2,s)}=3$ and $p_{(1,4),(1,4)}^{(1,4)}=0$, we obtain $s\neq4$. By $\Gamma_{2,s}\Gamma_{1,g-1}=\{\Gamma_{3,s-1}\}$ and $p_{(1,4),(4,1)}^{(2,3)}=0$, one has $s\neq5$. Since $\partial(x_{2,2},x_{0,0})>s-2$, we get $\partial(x_{4,2},x_{0,0})>s-4$, which implies $s=7$ from $\Gamma_{2,s}\Gamma_{2,3}=\{\Gamma_{4,s-2}\}$.

By Lemma \ref{jcs1_neq_2,s} (i), we get $\wz{\partial}(x_{0,0},x_{2,2})=(4,6)$. It follows from Lemma \ref{jiben} (i) and (iii) that $k_{4,6}=1$ and $\wz{\partial}(x_{0,0},x_{3,2})=(5,5)$. Since $x_{2,2}\in P_{(4,6),(1,4)}(x_{0,0},x_{3,2})$, we obtain $\partial(x_{0,4},x_{0,0})=1$. By Lemma \ref{jcs1_neq_2,s} (ii), $\Gamma_{1,4}\Gamma_{2,3}=\{\Gamma_{3,2},\Gamma_{3,6},\Gamma_{\wz{f}_{2}}\}$. Since $x_{0,4}\in P_{(2,7),(1,4)}(x_{4,3},x_{0,0})$ and $\Gamma_{2,7}\Gamma_{1,4}=\{\Gamma_{3,6}\}$, $(x_{4,2}',x_{0,0})\in\Gamma_{\wz{f}_{2}}$. The fact that $x_{1,1}\in\Gamma_{2,7}(x_{0,0})\cap\Gamma_{2,3}(x_{4,1})\cap\Gamma_{\wz{f}_{2}^{*}}(x_{3,1}')$, $x_{0,2}\in\Gamma_{2,7}(x_{4,1})\cap\Gamma_{2,3}(x_{0,0})$ and $x_{4,2}'\in\Gamma_{2,7}(x_{3,1}')\cap\Gamma_{\wz{f}_{2}^{*}}(x_{0,0})$ imply $\partial(x_{0,0},x_{4,1})=\partial(x_{4,1},x_{0,0})$ and $\partial(x_{0,0},x_{3,1}')=\partial(x_{3,1}',x_{0,0})$. Hence, $\partial(x_{3,1}',x_{0,0})=3$ or $\partial(x_{4,1},x_{0,0})=3$, contrary to $s=7$. Then $g\geq6$.

If $t=2$, from the claim and Lemma \ref{jcs1_neq_2,s} (iii), then $\partial(x_{0,0}',x_{3,0}')=\partial(x_{1,0},x_{3,0}')=3$, which implies $\partial(x_{0,0},x_{3,0}')=4$ since $\partial(x_{0,1},x_{3,0}')\geq3$; if $t\geq3$, then $\partial(x_{0,0},x_{3,0}')=4$. By $\wz{\partial}(x_{0,0},x_{4,0})=(4,g-4)$, one has $\partial(x_{0,0},x_{4,1})=5$. Hence, $\partial(x_{0,0},x_{4,2})=\partial(x_{0,0},x_{3,1}')+1$. Since $x_{2,1}\in P_{(2,s),(2,g-2)}(x_{1,0},x_{3,1}')$, from (i) and Lemma \ref{jiben} (iii), we get $\partial(x_{1,0},x_{3,1}')=4$. The fact that $\wz{\partial}(x_{0,0},x_{3,0}')=\wz{\partial}(x_{0,1},x_{3,1}')$ implies $\partial(x_{0,0},x_{4,2})=\partial(x_{0,0}',x_{3,1}')+2$. Thus, (ii) is valid.$\qed$

Pick a shortest path $(x_{3,3}=y_{0,j},y_{1,j},\ldots,y_{s_{j},j}=x_{2-j,2-j})$ for $0\leq j\leq r$, where $r=\max\{c\mid\partial(x_{0,0},x_{3+c,1+c})=4+2c~\textrm{and}~0\leq c\leq 2\}$. By Lemma \ref{jiben} (i), one has $k_{2+2j,s_{j}}=1$, which implies $|\Gamma_{2+2j,s_{j}}\Gamma_{t,g-t}|=1$ from Lemma \ref{jiben} (iii). In view of Lemma \ref{jcs1_neq_2,s} (i), we may assume that $\wz{\partial}(y_{0,0},y_{t,0})=(t,g-t)$ and $y_{i,0}=y_{i,1}=\cdots=y_{i,r}$ for $0\leq i\leq t$.

\begin{lemma}\label{main lemma}
Suppose $\partial(x_{0,0},x_{3+j,1+j})=4+2j$ for $j=1$ or $2$. If there exist two vertices $u_{1}$ and $y_{i,j'}$ such that $\partial(x_{3,3},u_{1})=1$ and $\partial(u_{1},y_{i,j'})=\partial(y_{1,j'},y_{i,j'})$ with $u_{1}\neq y_{1,j'}$, then $\wz{\partial}(x_{3,3},y_{i,j'})\neq\wz{\partial}(x_{3,3},y_{i,j''})$ for $0\leq j',j''\leq j$ and $j'\neq j''$.
\end{lemma}
\textbf{Proof.}~Suppose for the contrary that $\wz{\partial}(x_{3,3},y_{i,j'})=\wz{\partial}(x_{3,3},y_{i,j''})$. Since $k_{2+2j'',s_{j''}}=k_{2+2j',s_{j'}}=1$, from Lemma \ref{jiben} (iii), there exists a shortest path $(y_{0,j'},y_{1,j'},\ldots,y_{i,j'}=v_{i},v_{i+1},\ldots,v_{s_{j''}}=x_{2-j'',2-j''})$.

By the commutativity of $\Gamma$, there exist distinct vertices $u_{i-1}$ and $v_{i-1}$ such that $u_{i-1},v_{i-1}\in\Gamma_{g-1,1}(y_{i,j'})$ and $\partial(x_{3,3},u_{i-1})=\partial(x_{3,3},v_{i-1})=i-1$. It follows from Lemma \ref{jcs1_neq_2,s} (i) that $y_{i+1,j'},v_{i+1}\notin\Gamma_{2,s}(u_{i-1})\cup\Gamma_{2,s}(v_{i-1})$. Since $p_{(2,s),(g-1,1)}^{(1,g-1)}=1$, $y_{i+1,j'}=v_{i+1}$. By induction, we have $y_{c,j'}=v_{c}$ for $i\leq c\leq\min\{s_{j'},s_{j''}\}$. Without loss of generality, we may assume $s_{j'}<s_{j''}$. Since $k_{2+2j',s_{j'}}=1$, from Lemma \ref{jiben} (iii), we get $|\Gamma_{s_{j'},2+2j'}\Gamma_{g-1,1}|=1$ and $\partial(x_{3,3},v_{s_{j'}-1})=\partial(x_{3,3},b_{s_{j'}-1})$ with $b_{s_{j'}-1}\in P_{(g-1,1),(2,s)}(v_{s_{j'}},v_{s_{j'}+1})$, contrary to Lemma \ref{jcs1_neq_2,s} (i). The desired result holds.$\qed$

Write $\wz{d}_{t}:=\wz{\partial}(x_{0,1},x_{t+1,1}')$ and $\wz{h}_{t}:=\wz{\partial}(x_{0,0}',x_{t+1,1}')$.

\begin{lemma}\label{chain}
If $\partial(x_{0,0},x_{4,2})=6$, then $\partial(x_{0,0},x_{5,3})<8$, $\wz{d}_{t}\neq\wz{h}_{t}$, $\wz{\partial}(y_{0,j},y_{t+2,j})=\wz{h}_{t}$ and $\wz{\partial}(y_{i,j'},y_{i+t+2,j'})=\wz{d}_{t}$ for some $i\in\{0,1,\ldots, s_{j'}-t-2\}$ with $\{j,j'\}=\{0,1\}$.
\end{lemma}
\textbf{Proof.}~Let $r=\max\{c\mid\partial(x_{0,0},x_{3+c,1+c})=4+2c~\textrm{and}~1\leq c\leq2\}$.

For $0\leq \mu\leq r$, we claim $(y_{\lambda,\mu},y_{\lambda+t+2,\mu})\in\Gamma_{\wz{d}_{t}}$ for some $\lambda\in\{0,1,\ldots,s_{\mu}-t-2\}$ if $(y_{\lambda',\mu},y_{\lambda'+t+1,\mu})\notin\Gamma_{\wz{f}_{t}}$ for some $\lambda'\in\{0,1,\ldots,s_{\mu}-t-1\}$. It follows from Lemma \ref{jcs1_neq_2,s} (i) that $\wz{\partial}(y_{\lambda'',\mu},y_{\lambda''+t,\mu})=(t,g-t)$ for $0\leq\lambda''\leq s_{\mu}-t$. In view of $k_{2+\mu,s_{\mu}}=1$, Lemma \ref{jiben} (iii) and the commutativity of $\Gamma$, we may assume $\lambda'=0$. Since $\wz{\partial}(x_{0,0},x_{1,1})=(2,s)$ and $\partial(x_{1,1},x_{t,1})=t-1$, we have $(y_{0,\mu},y_{t+1,\mu})\in\Gamma_{t+1,g-t-1}$. By Lemma \ref{jcs1_neq_2,s} (iii), there exists an integer $\lambda$ such that $\lambda+1=\min\{c\mid\wz{\partial}(y_{c,\mu},y_{c+t+1,\mu})=\wz{f}_{t}\}$. Observe $0\leq \lambda\leq s_{\mu}-t-2$ and $(y_{\lambda,\mu},y_{\lambda+t+1,\mu})\in\Gamma_{t+1,g-t-1}$. Since $x_{1,2}\in P_{(2,s),(t,g-t)}(x_{0,1},x_{t+1,2})$, from Lemma \ref{jcs1_neq_2,s} (i), we get $\wz{\partial}(y_{\lambda,\mu},y_{\lambda+t+2,\mu})=\wz{d}_{t}$. Thus, the claim is valid.

Let $0\leq \mu,\mu'\leq r$ and $\mu\neq \mu'$. Suppose $(y_{\lambda,\mu},y_{\lambda+t+1,\mu}),(y_{\lambda',\mu'},y_{\lambda'+t+1,\mu'})\in\Gamma_{\wz{f}_{t}}$ for $0\leq\lambda\leq s_{\mu}-t-1$ and $0\leq\lambda'\leq s_{\mu'}-t-1$. By Lemma \ref{jcs1_neq_2,s} (i), one has $(y_{\lambda,\mu},y_{\lambda+t,\mu}),(y_{\lambda',\mu'},y_{\lambda'+t,\mu'})\in\Gamma_{t,g-t}$ for $0\leq \lambda\leq s_{\mu}-t$ and $0\leq \lambda'\leq s_{\mu'}-t$. In view of Lemma \ref{jiben} (iv) and Lemma \ref{jcs1_neq_2,s} (ii), we get $p_{\wz{f}_{t},(g-1,1)}^{(t,g-t)}=1$, which implies that $y_{a,\mu}=y_{a,\mu'}$ for $0\leq a\leq\min\{s_{\mu},s_{\mu'}\}$.  Without loss of generality, we may assume $s_{\mu}<s_{\mu'}$. Since $k_{2+2\mu,s_{\mu}}=1$, from Lemma \ref{jiben} (iii), we obtain $|\Gamma_{s_{\mu},2+2\mu}\Gamma_{g-1,1}|=1$ and $\partial(x_{3,3},y_{s_{\mu}-1,\mu'})=\partial(x_{3,3},b_{s_{\mu}-1})$ with $b_{s_{\mu}-1}\in P_{(g-1,1),(2,s)}(y_{s_{\mu},\mu'},y_{s_{\mu}+1,\mu'})$, contrary to Lemma \ref{jcs1_neq_2,s} (i).

Suppose that $(y_{\lambda,\mu},y_{\lambda+t+1,\mu})\notin\Gamma_{\wz{f}_{t}}$ and $(y_{\lambda',\mu'},y_{\lambda'+t+1,\mu})\notin\Gamma_{\wz{f}_{t}}$ for some $\lambda\in\{0,1,\ldots,s_{\mu}-t-1\}$ and $\lambda'\in\{0,1,\ldots,s_{\mu'}-t-1\}$. In view of the claim, Lemma \ref{jiben} (iii) and the commutativity of $\Gamma$, we may assume $y_{t+2,\mu},y_{t+2,\mu'}\in\Gamma_{\wz{d}_{t}}(x_{3,3})$. Observe $x_{t+1,1}\in P_{(t+1,g-t-1),(1,g-1)}(x_{0,1},x_{t+1,1}')$ and $x_{1,1}\in P_{(1,g-1),\wz{f}_{t}}(x_{0,1},x_{t+1,1}')$, contrary to Lemma \ref{main lemma}. By the claim again, we get $r=1$ and $\partial(x_{0,0},x_{5,3})<8$. Then $(y_{i,j},y_{i+t+1,j})\in\Gamma_{\wz{f}_{t}}$ for $0\leq i\leq s_{j}-t-1$, and $(y_{i',j'},y_{i'+t+2,j'})\in\Gamma_{\wz{d}_{t}}$ for some $i'\in\{0,1,\ldots, s_{j'}-t-2\}$  with $\{j,j'\}=\{0,1\}$.

By Lemma \ref{jcs1_neq_2,s} (iii) and Lemma \ref{main lemma}, one has $p_{(1,g-1),(t,g-t)}^{\wz{f}_{t}}=1$ and $|\Gamma_{\wz{f}_{t}}\Gamma_{g-1,1}|\neq1$. In view of Lemma \ref{jiben} (iii), we obtain $k_{\wz{f}_{t}}\neq1$, which implies $t+1<s_{j}$ by $k_{2+2j,s_{j}}=1$. Since $(x_{1,1},x_{t+1,1}')\in\Gamma_{\wz{f}_{t}}$ and $x_{2,1}\in P_{(2,s),(t,g-t)}(x_{1,0},x_{t+1,1}')$, from Lemma \ref{jcs1_neq_2,s} (i), we get $\wz{\partial}(y_{0,j},y_{t+2,j})\neq\wz{\partial}(x_{1,0},x_{t+1,1}')$. Observe $x_{1,1}\in P_{(1,g-1),\wz{f}_{t}}(x_{0,1},x_{t+1,1}')$ and $x_{t+1,1}\in P_{(t+1,g-t-1),(1,g-1)}(x_{0,1},x_{t+1,1}')$. By Lemma \ref{main lemma}, Lemma \ref{jiben} (iii) and the commutativity of $\Gamma$, we have $\wz{\partial}(y_{0,j},y_{t+2,j})=\wz{h}_{t}$ and $\wz{h}_{t}\neq\wz{d}_{t}$.$\qed$

\noindent\textbf{Proof of Lemma \ref{pre-main}.}

\textbf{Case 1.} $t=1$.

By Lemma \ref{jcs1_neq_2,s} (ii), $(A_{1,g-1})^{2}=A_{2,g-2}+A_{2,l}+3A_{2,s}$ with $g-2<l<s$. Lemma \ref{jcs1d=4} (i) implies $g\geq5$. It follows from the claim in Lemma \ref{jcs1d=4} that $\wz{\partial}(x_{0,2}',x_{2,2}')=\wz{\partial}(x_{1,1}',x_{3,1}')=(2,g-2)$ and there exist two paths $(x_{0,2}',z_{0},x_{2,2}')$ and $(x_{1,1}',z_{1},x_{3,1}')$. Since $p_{(1,g-1),(1,g-1)}^{(2,s)}=3$, we get $z_{0}\in P_{(1,g-1),(1,g-1)}(x_{1,1}',x_{2,2}')$. By $p_{(g-1,1),(2,g-2)}^{(1,g-1)}=1$, we have $\wz{\partial}(x_{3,1}',x_{5,2})\neq(2,g-2)$, which implies $\partial(x_{0,2}',x_{3,1}')>2$ and $z_{0}\neq z_{1}$. Since $p_{(1,g-1),(1,g-1)}^{(2,g-2)}=1$, we obtain $\Gamma_{1,g-1}(x_{1,1}')=\{z_{0},z_{1},x_{2,2}\}$.

Since $x_{4,1}\in P_{(3,g-3),(1,g-1)}(x_{1,1},x_{4,1}')$, from the commutativity of $\Gamma$, we have $x_{1,2}$ or $x_{1,1}'\in P_{(1,g-1),(3,g-3)}(x_{1,1},x_{4,1}')$. The fact that $p_{(g-1,1),(2,g-2)}^{(1,g-1)}=1$ implies $(x_{4,1}',x_{6,2})\notin\Gamma_{2,g-2}$ and $x_{1,2}\notin\Gamma_{g-3,3}(x_{4,1}')$. Then $x_{1,1}'\in\Gamma_{g-3,3}(x_{4,1}')$.

Pick a path $(x_{1,1}',w_{0},w_{1},x_{4,1}')$. Since $(x_{4,1}',x_{6,2})\notin\Gamma_{2,g-2}$, $w_{0}\in\{z_{0},z_{1}\}$. Observe $(x_{1,1}',x_{3,1}')\in\Gamma_{2,g-2}$ and $(x_{1,1},x_{3,1}')\notin\Gamma_{3,g-3}$. By $p_{(2,g-2),(g-1,1)}^{(1,g-1)}=1$, we have $(x_{1,1},z_{1})\in\Gamma_{2,l}$, which implies $(x_{1,1},z_{0})\in\Gamma_{2,g-2}$. Since $(x_{1,1},x_{4,1}')\notin\Gamma_{4,g-4}$, one gets $w_{0}=z_{1}$. By $p_{(2,g-2),(g-1,1)}^{(1,g-1)}=1$ again, we obtain $w_{1}=x_{3,1}'$, which implies $x_{4,1},x_{3,1}'\in P_{(1,g-1),(1,g-1)}(x_{3,1},x_{4,1}')$, contrary to $p_{(1,g-1),(1,g-1)}^{(2,l)}=1$.

\textbf{Case 2.} $t=2$.

We claim $|\Gamma_{\wz{f}_{2}}\Gamma_{1,g-1}|=3$ and $\Gamma_{\wz{h}_{2}}\notin\Gamma_{3,g-3}\Gamma_{1,g-1}$. By Lemma \ref{jcs1_neq_2,s} (i), Lemma \ref{jcs1d=4} (i) and Lemma \ref{jiben} (iii), we get $\Gamma_{2,s}\Gamma_{2,g-2}=\{\Gamma_{4,s-2}\}$. Then $(x_{0,1},x_{3,2}),(x_{1,0},x_{3,1}')\in\Gamma_{4,s-2}$. Lemma \ref{jcs1d=4} (ii) implies $g\geq6$ and $p_{(2,g-2),(1,g-1)}^{(3,g-3)}=p_{(g-2,2),(g-1,1)}^{(g-3,3)}=1$.  Hence, $\wz{d}_{2}\neq(4,s-2)$ and $\Gamma_{3,g-3}\Gamma_{1,g-1}=\{\Gamma_{4,g-4},\Gamma_{4,s-2},\Gamma_{\wz{d}_{2}}\}$. Since $(x_{1,1},x_{3,1}')\in\Gamma_{\wz{f}_{2}}$, we obtain $\Gamma_{4,s-2},\Gamma_{\wz{d}_{2}}\in\Gamma_{\wz{f}_{2}}\Gamma_{1,g-1}$. If $\partial(x_{0,0},x_{4,2})=6$, by Lemma \ref{jcs1_neq_2,s} (i) and Lemma \ref{chain}, then $\wz{h}_{2}\notin\{(4,g-4),(4,s-2),\wz{d}_{2}\}$; if $\partial(x_{0,0},x_{4,2})<6$, by Lemma \ref{jcs1d=4} (ii), then $\partial(x_{0,1},x_{3,1}')=\partial(x_{0,0},x_{3,0}')=4$ and $\partial(x_{0,0}',x_{3,1}')<4$. The claim is valid.

In view of Lemma \ref{jiben} (iv), we have $(A_{1,g-1})^{2}=A_{2,g-2}+3A_{2,s}$, which implies $k_{2,g-2}=6$ from Lemma \ref{jiben} (i) and (vi). Since $p_{(2,g-2),(1,g-1)}^{(3,g-3)}=1$, by Lemma \ref{jcs1_neq_2,s} (ii), we get $k_{3,g-3}=6$. In view of $\Gamma_{2,s}\Gamma_{1,g-1}=\{\Gamma_{3,s-1}\}$, one obtains $k_{3,s-1}=3$.

By the claim in Lemma \ref{jcs1d=4}, one has $\wz{\partial}(x_{0,1}',x_{3,1}')=(3,g-3)$ and there exists a path $(x_{0,1}',u_{1},u_{2},x_{3,1}')$. From $p_{(1,g-1),(1,g-1)}^{(2,s)}=3$, $u_{2}\in P_{(1,g-1),(1,g-1)}(x_{2,0}',x_{3,1}')$. Since $p_{(1,g-1),(2,g-2)}^{(3,g-3)}=1$, one gets $u_{1}\neq x_{1,2}$ and $\wz{\partial}(x_{0,1},u_{1})=(2,g-2)$. By $(x_{0,1},x_{3,1}')\notin\Gamma_{4,g-4}$, $p_{(3,g-3),(g-1,1)}^{(2,g-2)}=1$ and $\Gamma_{2,s}\Gamma_{1,g-1}=\{\Gamma_{3,s-1}\}$, we obtain $\wz{\partial}(x_{0,1},u_{2})=\wz{\partial}(x_{1,1},x_{3,1}')=\wz{f}_{2}$. Since $\Gamma_{2,s}\Gamma_{2,g-2}=\{\Gamma_{4,s-2}\}$, we have $(x_{g-1,0}',u_{2})\in\Gamma_{4,s-2}$, which implies $\wz{\partial}(x_{-1,1},u_{2})=\wz{h}_{2}$ or $\wz{d}_{2}$ from $(x_{1,0},x_{3,1}')\in\Gamma_{4,s-2}$.

\textbf{Case 2.1.} $\partial(x_{0,0},x_{4,2})<6$.

In view of Lemma \ref{jcs1d=4} (ii), $\partial(x_{0,0}',x_{3,1}')<4$. Since $\Gamma_{2,s}\Gamma_{2,g-2}=\{\Gamma_{4,s-2}\}$, we get $(x_{1,1}',x_{3,1}')\notin\Gamma_{2,g-2}$. By $p_{(2,g-2),(1,g-1)}^{(3,g-2)}=1$, one has $(x_{1,2},x_{3,1}')\notin\Gamma_{2,g-2}$. Hence, $p_{(1,g-1),(2,g-2)}^{\wz{f}_{2}}=1$. Lemma \ref{jcs1_neq_2,s} (ii) and Lemma \ref{jiben} (vi) imply $k_{\wz{f}_{2}}=6$.

Since $k_{2,s}=1$, from Lemma \ref{jiben} (iii) and Lemma \ref{jcs1_neq_2,s} (iii), we may assume $\wz{\partial}(y_{0,0},y_{3,0})=\wz{\partial}(x_{1,1},x_{3,1}')$. By Lemma \ref{jcs1_neq_2,s} (i) and $\partial(x_{0,0}',x_{3,1}')<4$, $\wz{\partial}(y_{0,0},y_{4,0})=\wz{d}_{2}$. Observe $x_{1,1}\in P_{(2,s),\wz{f}_{2}}(x_{0,0},x_{3,1}')$ and $x_{0,1}'\in P_{(2,s),(3,g-3)}(x_{g-1,0}',x_{3,1}')$. Since $s>4$, we get $\wz{\partial}(y_{0,0},y_{5,0})=\wz{\partial}(x_{-1,1},x_{3,1}')$, which implies $\wz{\partial}(x_{-1,1},u_{2})=\wz{d}_{2}$ by $\partial(x_{0,0}',x_{3,1}')<4$. Hence, $\wz{\partial}(x_{0,0},u_{2})=\wz{h}_{2}$.

Since $\Gamma_{2,s}\Gamma_{1,g-1}=\{\Gamma_{3,s-1}\}$ and $p_{(1,g-1),(1,g-1)}^{(2,g-2)}=1$, from the claim, we have $\wz{\partial}(x_{0,0},u_{1})=\wz{f}_{2}$, which implies $u_{1},x_{2,0}'\in P_{\wz{f}_{2},(1,g-1)}(x_{0,0},u_{2})$. In view of Lemma \ref{jcs1_neq_2,s} (ii), one gets $(\Gamma_{1,g-1})^{3}=\{\Gamma_{3,g-3},\Gamma_{3,s-1},\Gamma_{\wz{f}_{2}}\}$. Since $\partial(x_{0,0}',x_{3,1}')<4$, from the claim and Lemma \ref{jiben} (vi), we obtain $p_{\wz{f}_{2},(1,g-1)}^{\wz{h}_{2}}=2$ and $k_{\wz{h}_{2}}=3$. Hence, $\wz{h}_{2}=(1,g-1)$ or $(3,s-1)$. If $\wz{h}_{2}=(1,g-1)$, by Lemma \ref{jiben} (i) and Lemma \ref{jcs1_neq_2,s} (iii), then $\wz{f}_{2}=(3,3)$ since $k_{\wz{f}_{2}}=6$ and $x_{0,0}'\in P_{(g-1,1),(1,g-1)}(x_{1,1},x_{3,1}')$, which implies $g<6$, a contradiction. If $\wz{h}_{2}=(3,s-1)$, from $\Gamma_{3,s-1}\in\Gamma_{2,s}\Gamma_{1,g-1}$, then $\partial(x_{1,1}',x_{3,1}')=1$, contrary to $\partial(y_{0,0},y_{3,0})=\partial(x_{1,1},x_{3,1}')=3$.

\textbf{Case 2.2.} $\partial(x_{0,0},x_{4,2})=6$.

In view of Lemma \ref{jcs1_neq_2,s} (iii) and Lemma \ref{main lemma}, we obtain $p_{(1,g-1),(2,g-2)}^{\wz{f}_{2}}=1$. It follows from Lemma \ref{jcs1_neq_2,s} (ii) and Lemma \ref{jiben} (vi) that $k_{\wz{f}_{2}}=6$. Since $g\geq6$, we have $p_{(3,g-3),(1,g-1)}^{(4,g-4)}=p_{(g-3,3),(g-1,1)}^{(g-4,4)}=1$, which implies $\partial(x_{0,2},x_{3,1}')>3$ and $p_{(1,g-1),\wz{f}_{2}}^{\wz{d}_{2}}=1$. In view of the claim and Lemma \ref{jiben} (vi) again, one gets $k_{\wz{d}_{2}}=6$.

Note that $k_{2,s}=k_{4,s_{1}}=1$. By Lemma \ref{chain}, Lemma \ref{jiben} (iii) and the commutativity of $\Gamma$, we may assume $(y_{0,j},y_{4,j})\in\Gamma_{\wz{h}_{2}}$ and $(y_{0,j'},y_{4,j'})\in\Gamma_{\wz{d}_{2}}$ with $\{j,j'\}=\{0,1\}$. Observe $x_{1,1}\in P_{(2,s),\wz{f}_{2}}(x_{0,0},x_{3,1}')$ and $x_{0,1}'\in P_{(2,s),(3,g-3)}(x_{g-1,0}',x_{3,1}')$. By Lemma \ref{jcs1_neq_2,s} (i), we obtain $\wz{\partial}(y_{0,j'},y_{5,j'})=\wz{\partial}(x_{-1,1},x_{3,1}')$. Since $k_{2+2j',s_{j}'}=1$ and $k_{\wz{d}_{2}}=6$, from Lemma \ref{jiben} (vi), one has $s_{j}'>5$. Observe $x_{0,1}\in P_{(1,g-1),\wz{d}_{2}}(x_{-1,1},x_{3,1}')$ and $x_{3,1}\in P_{(4,g-4),(1,g-1)}(x_{-1,1},x_{3,1}')$. Since $p_{(2,s),(g-1,1)}^{(1,g-1)}=1$, by Lemma \ref{chain} and Lemma \ref{jcs1_neq_2,s} (i), we get $(x_{-1,1},u_{2})\in\Gamma_{\wz{d}_{2}}$ and $(x_{0,0},u_{2})\in\Gamma_{\wz{h}_{2}}$.

Since $\Gamma_{2,s}\Gamma_{1,g-1}=\{\Gamma_{3,s-1}\}$, from the claim, we get $u_{1},x_{2,0}'\in P_{\wz{f}_{2},(1,g-1)}(x_{0,0},u_{2})$. In view of $k_{2+2j,s_{j}}=1$, $k_{\wz{f}_{2}}=6$ and Lemma \ref{jiben} (vi), one has $s_{j}>4$. By $p_{(2,s),(g-1,1)}^{(1,g-1)}=1$ and Lemma \ref{jcs1_neq_2,s} (i), we obtain $p_{(1,g-1),\wz{f}_{2}}^{\wz{h}_{2}}=2$ and $\wz{\partial}(y_{1,j},y_{4,j})=\wz{\partial}(y_{0,j},y_{3,j})=\wz{f}_{2}$. It follows from Lemma \ref{jiben} (vi) that $k_{\wz{h}_{2}}=3$. Pick a vertex $v_{1}\in P_{(1,g-1),\wz{f}_{2}}(y_{0,j},y_{4,j})\setminus\{y_{1,j}\}$ and a path $(v_{1},v_{2},v_{3},y_{4,j})$. Observe $x_{1,1}\in P_{(1,g-1),\wz{f}_{2}}(x_{0,1},x_{3,1}')$ and $x_{3,1}\in P_{(3,g-3),(1,g-1)}(x_{0,1},x_{3,1}')$. By Lemma \ref{main lemma}, Lemma \ref{jiben} (iii) and the commutativity of $\Gamma$, we have $y_{1,j},v_{1}\notin\Gamma_{\wz{d}_{2}^{*}}(y_{5,j})$. Since $x_{2,1}\in P_{(2,s),(2,g-2)}(x_{1,0},x_{3,1}')$ and $(x_{1,1},x_{3,1}')\in\Gamma_{\wz{f}_{2}}$, from Lemma \ref{jcs1_neq_2,s} (i), one gets $v_{1},y_{1,j}\in P_{(1,g-1),\wz{h}_{2}}(y_{0,j},y_{5,j})$. The fact $(y_{0,j},y_{3,j})\in\Gamma_{\wz{f}_{2}}$ and $p_{(1,g-1),(2,g-2)}^{\wz{f}_{2}}=1$ implies $v_{3}\neq y_{3,j}$. By $p_{(2,s),(g-1,1)}^{(1,g-1)}=1$ and Lemma \ref{jcs1_neq_2,s} (i), one obtains $p_{\wz{\partial}(y_{0,j},y_{5,j}),(g-1,1)}^{\wz{h}_{2}}$=1. From Lemma \ref{jiben} (ii), one has $k_{\wz{\partial}(y_{0,j},y_{5,j})}=1$ and $p_{(1,g-1),\wz{h}_{2}}^{\wz{\partial}(y_{0,j},y_{5,j})}=3$. Thus, $s_{j}=5$.

Since $k_{4,s_{1}}=1$, from Lemma \ref{jiben} (iii), one gets $\wz{\partial}(x_{0,0},x_{3,2})=(5,s_{1}-1)$. By $x_{2,2}\in P_{(4,s_{1}),(1,g-1)}(x_{0,0},x_{3,2})$ and Lemma \ref{jiben} (i), we have $j=0$ and $g=6$. Since $\Gamma_{2,5}\Gamma_{2,4}=\{\Gamma_{4,3}\}$, we obtain $k_{4,3}=6$, contrary to $k_{3,4}=k_{3,s-1}=3$.

\textbf{Case 3.} $t\geq3$.

Observe $x_{1,1}\in\Gamma_{2,s}(x_{0,0})\cap\Gamma_{g-3,3}(x_{4,1})\cap\Gamma_{g-3,3}(x_{3,1}')$. By Lemma \ref{jiben} (iii) and Lemma \ref{jcs1d=4} (ii), we have $\partial(x_{0,0},x_{3,1}')=5$, which implies $\partial(x_{0,0},x_{4,2})=6$ and $\partial(x_{0,0},x_{5,3})<8$ from Lemma \ref{chain}.

Since $k_{4,s_{1}}=1$ and $\partial(x_{0,0},x_{3,2})=5$, from Lemma \ref{jiben} (iii), we get $\partial(x_{0,0},x_{3,3})=6$. By $x_{4,2},x_{3,2}'\in\Gamma_{2,g-2}(x_{2,2})$, one has $\partial(x_{0,0},x_{3,2}')=6$, which implies $\partial(x_{0,0},x_{4,3})=7$. Observe that $x_{2,2}\in P_{(4,s_{1}),(3,g-3)}(x_{0,0},x_{5,2})$ and $x_{2,2}\in P_{(4,s_{1}),(3,g-3)}(x_{0,0},x_{4,2}')$. By Lemma \ref{jiben} (iii) and $\partial(x_{0,0},x_{5,3})<8$, we obtain $\partial(x_{0,0},x_{5,2})<7$.

If $t\geq5$, then $(x_{0,0},x_{4,0}')\in\Gamma_{5,g-5}$; if $t=4$, then $(x_{0,0},x_{4,0}')\in\Gamma_{\wz{f}_{4}}$; if $t=3$, then $(x_{0,0},x_{4,0}')\in\Gamma_{\wz{d}_{3}}$. By Lemma \ref{jcs1_neq_2,s} (iii) or Lemma \ref{chain}, we have $\partial(x_{0,0},x_{4,0}')=5$. Since $\partial(x_{0,0},x_{5,0})=\partial(x_{0,0},x_{4,1})=5$, we have $\partial(x_{0,0},x_{5,1})=6$. Hence, $\partial(x_{0,0},x_{4,1}')<6$.

Since $x_{2,1}\in P_{(2,s),(3,g-3)}(x_{1,0},x_{4,1}')$, from Lemma \ref{jiben} (iii), we get $\partial(x_{1,0},x_{4,1}')=\partial(x_{0,0},x_{4,1})=5$. By $\wz{\partial}(x_{0,0},x_{4,0}')=\wz{\partial}(x_{0,1},x_{4,1}')$, one has $\partial(x_{0,0}',x_{4,1}')<5$, which implies $t=3$ or $4$, contrary to Lemma \ref{chain} or Lemma \ref{jcs1_neq_2,s} (iii).$\qed$

\section*{Acknowledgement}
The authors would like to thank the anonymous reviewers for their careful reading of the manuscript of the paper and their invaluable, critical detailed suggestions which led to a great improvement of the presentation of the paper. Y. Yang is supported by the Fundamental Research Funds for the Central Universities, B. Lv is supported by NSFC (11501036), K. Wang is supported by NSFC (11671043, 11371204), and the Fundamental Research Funds for the Central Universities.

\begin{table}[hbt]
\begin{center}
\renewcommand\arraystretch{1.15}
\caption{\quad Two way distance of digraphs in Theorem \ref{Main}} \vspace{0.2cm}

\setlength{\tabcolsep}{1.8pt}
\begin{tabular}{|c|c|c|}
\hline
$\Gamma$ & Conditions & $\wz{\partial}((0,0),(a,b))$ with $(a,b)\neq(0,0)$ \\\hline
\multirow{2}{*}{(vi)} & $a=0$ & $(g,g)$ \\\cline{2-3}
& $a\neq0$ & $(\hat{a},g-\hat{a})$ \\\hline
\multirow{5}{*}{(vii)} & $b=0$ & $(\min\{\hat{a},2n-2\hat{a}\},\min\{n-\hat{a},2\hat{a}\})$ \\\cline{2-3}
& $a=0$ & $(\min\{\hat{b},2n-2\hat{b}\},\min\{n-\hat{b},2\hat{b}\})$ \\\cline{2-3}
& $a=b$ & $(\min\{n-\hat{a},2\hat{a}\},\min\{\hat{a},2n-2\hat{a}\})$ \\\cline{2-3}
& $0<\hat{b}<\hat{a}$ & $(h_{0},l_{0})$ \\\cline{2-3}
& $0<\hat{a}<\hat{b}$ & $(h_{1},l_{1})$ \\\hline
\multirow{4}{*}{(viii)} & $-n<\hat{b}-\hat{a}<0$ & $(\min\{3n-3\hat{a}+\hat{b},3\hat{a}-2\hat{b}\},\min\{3\hat{a}-\hat{b},3n+2\hat{b}-3\hat{a}\})$ \\\cline{2-3}
& $0\leq\hat{b}-\hat{a}<n$ & $(\hat{b},\max\{3\hat{a}-\hat{b},2\hat{b}-3\hat{a}\})$ \\\cline{2-3}
& $n\leq\hat{b}-\hat{a}<2n$ & $(\max\{3n+3\hat{a}-2\hat{b},\hat{b}-3\hat{a}\},3n-\hat{b})$ \\\cline{2-3}
& $2n\leq\hat{b}-\hat{a}<3n$ & $(\min\{6n+3\hat{a}-2\hat{b},\hat{b}-3\hat{a}\},\min\{3n+3\hat{a}-\hat{b},2\hat{b}-3\hat{a}-3n\})$ \\\hline
\end{tabular}
\end{center}

For any element $a$ in a residue class ring, we assume that $\hat{a}$ denotes the minimum nonnegative integer in $a$. For $i=0,1$, let
\begin{center}
$h_{i}=\min\{\hat{a}+\hat{b},(i+1)n+\hat{a}-2\hat{b},(2-i)n-2\hat{a}+\hat{b}\}$,\\
$l_{i}=\min\{2n-\hat{a}-\hat{b},(1-i)n-\hat{a}+2\hat{b},in-\hat{b}+2\hat{a}\}$.
\end{center}
\end{table}

\end{document}